\input epsf

 \documentclass[12pt,a4paper]{article}
 \usepackage[dvips]{graphicx}
 \usepackage{latexsym,amsmath,amsfonts,amscd,subfigure}
 \usepackage{epsfig}
 \usepackage{changebar}
 \usepackage{indentfirst}
 \usepackage{verbatim}
 \usepackage{color}
 \usepackage{colordvi}

\begin{document}
\baselineskip=2pc

\begin{center}
{\Large \bf  Moment-based multi-resolution HWENO scheme for hyperbolic conservation laws}
\end{center}

\centerline{Jiayin Li\footnote{School of Mathematical Sciences, Xiamen University, Xiamen, Fujian 361005, P.R. China. E-mail: jiayin@stu.xmu.edu.cn.},
Chi-Wang Shu\footnote{Division of Applied Mathematics, Brown University, Providence, RI 02912. E-mail:
Chi-Wang\_Shu@brown.edu.}
and Jianxian Qiu\footnote{School of Mathematical Sciences and Fujian Provincial Key Laboratory of Mathematical Modeling and High-Performance
Scientific Computing, Xiamen University, Xiamen, Fujian 361005, P.R. China. E-mail: jxqiu@xmu.edu.cn.
}
}

\vspace{.1in}
\centerline{\bf Abstract}
\baselineskip=1.7pc

In this paper, a high-order moment-based multi-resolution Hermite weighted essentially
non-oscillatory (HWENO) scheme is designed for hyperbolic conservation laws.
The main idea of this scheme is derived from our previous work [J. Comput. Phys., 446 (2021) 110653], in which the integral averages of the function and its first order derivative are used to reconstruct both the function and its first order derivative values at the boundaries.
However, in this paper, only the function values at the Gauss-Lobatto points in the one or two dimensional case need to be reconstructed by using the information of the zeroth and first order moments. In addition, an extra modification procedure is used
to modify those first order moments in the troubled-cells, which leads to an improvement of stability and an enhancement of resolution near discontinuities.
To obtain the same order of accuracy, the size of the stencil required by this moment-based multi-resolution HWENO scheme is still the same as the general HWENO scheme and is more compact than the general WENO scheme.
Moreover, the linear weights can also be any positive numbers as long as their sum equals one and the CFL number can still be 0.6 whether for the one or two dimensional case.
Extensive numerical examples are given to demonstrate the stability and resolution of such moment-based
multi-resolution HWENO scheme.

\vfill
{\bf AMS(MOS) subject classification:} 65M60, 35L65

 \pagenumbering{arabic}
 \newpage
 \baselineskip=2pc

\section{Introduction}
\label{sec1}
\setcounter{equation}{0}
\setcounter{figure}{0}
\setcounter{table}{0}
{
In this paper, a high-order moment-based multi-resolution Hermite weighted essentially non-oscillatory
(HWENO) scheme is designed for hyperbolic conservation laws
\begin{align}
\label{hyperbolic-conservation-laws}
\left\{
\begin{aligned}
&u_t+\nabla \cdot f(u)=0,\\
&u(x_1,...,x_d,0)=u_0(x_1,...,x_d).
\end{aligned}
\right.
\end{align}
We concentrate our attention on the one and two dimensional cases ($d=1$ or 2), and in these cases we denote
$x_1$ as $x$ and $x_2$ as $y$.

Conservation laws arise from the physical law that the conservative quantity in any control body can change
only due to the flux passing through its boundaries, which naturally hold for many fundamental physical quantities,
such as the mass, momentum, energy and so on.
Such conservation laws are widely used in a broad spectrum of disciplines where wave motion or advective
transport is important: gas dynamics, acoustics, elastodynamics, optics, geophysics, and biomechanics,
to name but a few.

The differential equation (\ref{hyperbolic-conservation-laws}) can be derived from the integral equation
by simple manipulations provided that the conservative quantity and its corresponding flux are
sufficiently smooth.
This proviso is important because in practice many interesting solutions are not smooth, but contain
discontinuities such as shock waves.
A fundamental feature of nonlinear conservation laws is that discontinuities can easily develop
spontaneously even from smooth initial data, and must be dealt with carefully both mathematically
and computationally.
At a discontinuity in the conservative quantity, the differential equation does not hold in the classical
sense and it is important to remember that the integral form of the conservation laws does continue to hold
which is more fundamental.
This is also why we choose conservative schemes, such as the finite volume method considered in this
paper, which is based on the integral form of the conservation laws.

Since conservation laws have a very wide range of applications and it is almost impossible in general to
get their exact solutions, many scholars have explored and proposed a series of numerical methods and are still
trying to improve the performance of these algorithms.
In 1994, Liu et al. proposed the first finite volume WENO scheme in \cite{loc1994}, and then, in 1996, Jiang and Shu
improved this WENO scheme to fifth order and to conservative finite difference formulation (which is more
efficient in multi-dimensions), and gave a general definition of the smoothness indicators
and nonlinear weights in \cite{js1996}.
The methodology of such WENO schemes is to use a nonlinear convex combination of all the candidate
stencils to improve the order of accuracy in smooth regions without destroying the non-oscillatory
behavior near discontinuities.
This is also the difference of such WENO schemes from the ENO schemes in \cite{heoc1987,so1988,so1989},
which only choose the locally smoothest stencil automatically among all the central and biased
spatial stencils.
Thereafter, different kinds of WENO schemes have been developed in, e.g.
\cite{bs2000,p2002,ccd2011,dbsr2017,zs2018,zs2019,bgs2016,bgfb2020}.
Although these WENO schemes work well for most of the problems we encountered, there is still room for
improvement. For example, if we want to obtain a higher order scheme, we must further expand the
stencil. This will make our scheme not very compact and will also bring trouble to the processing of
the boundary conditions.
In order to overcome this drawback, Qiu and Shu proposed the first HWENO scheme
for one-dimensional problems in \cite{qs2004} and then, in 2005, they extended
this HWENO scheme to two-dimensional problems in \cite{qs2005}, where two different stencils were used
to reconstruct the function and its first order derivative values, respectively.
The main difference of such HWENO scheme from the WENO scheme is that both the function and its first
order derivative values are evolved in time and used in the reconstruction process, not like the WENO
scheme in which only the function values are evolved and used. This allows the HWENO scheme to obtain
the same order of accuracy as the WENO scheme with relatively narrower stencils.
But there occurs a new issue, that is this HWENO scheme is not stable enough when simulating certain
severe problems with strong discontinuities, including the double Mach and forward step problems.
This difficulty is largely due to the fact that the first order derivative values may become very large
near these discontinuities.
Thus, the stability issue may arise, if these large values are used straightforwardly without any
modification.
Driven by the goal of solving this issue, many effective methods based on the idea of the original
HWENO scheme have emerged. For example,
the scheme with a new procedure to reconstruct the first order derivative values by Zhu and Qiu
in \cite{zq2008} in 2008,
the scheme with an additional positivity-preserving limiter by Liu and Qiu in \cite{lq2015,lq2016} in
2015 and 2016,
the scheme with a troubled-cell indicator to modify the first order moments near the discontinuities
before the reconstruction algorithm by Zhao et al. in \cite{zcq2020,zq2020} in 2020,
the scheme with a hierarchy of nested central spatial stencils by Li et al. in \cite{lqs2021} in 2021
and so on, have been developed.

In 2018 and 2019, Zhu and Shu proposed a new type of high-order finite difference and finite volume
multi-resolution WENO schemes in \cite{zs2018,zs2019}, which only use the information on a
hierarchy of nested central spatial stencils and can not only obtain high order accuracy in
smooth regions but also allow the accuracy near discontinuities to degrade gradually.
Following the idea of these multi-resolution WENO schemes, we designed a new type of high-order finite
volume and finite difference multi-resolution HWENO schemes in \cite{lqs2021} in 2021, for which only
the function values need to be reconstructed by the HWENO schemes, and the first order derivative
values are obtained directly from the polynomial with the highest degree in the hierarchy. This can
improve the resolution of the scheme but does not have much effect on its stability, since the first
two layers in the hierarchy do not contain the information of the first order derivative.  In this
paper, the function and its first order moment values are used in our reconstruction algorithm, and
only the function values are needed to be reconstructed according to our control equations. Thus, the
coupling between the function and its first order moment is stronger.

There are two issues to be addressed.
The first one is that,
although the bigger the linear weights are for higher degree polynomials, the steeper the shock transitions
become near the discontinuities, the gap between these linear weights cannot be too large, otherwise the
corresponding nonlinear weights will still be too close to the linear weights near the discontinuities.
This will cause problems, since the higher degree polynomials, which require the information of the
first order derivative or moment, account for too much in the final reconstruction polynomial, but
such first order derivative or moment values may be very large near the discontinuities.
The other one is that
with the suitable choice of the linear weights, the order of the final reconstruction polynomial will
degrade gradually near the discontinuities until it drops to first order, which will smooth out these
discontinuities to a certain extent.
Guided by the idea of Zhao et al. in \cite{zcq2020,zq2020}, we first perform the reconstruction algorithm, and then modify the first order moments of the troubled-cells and repeat the reconstruction algorithm for these troubled-cells to update the
corresponding Gauss-Lobatto point values.
After such a modification procedure, the proportion of the last two layers will become much higher than
that of the first two layers in the reconstruction process, thus the resolution near the
discontinuities can be increased significantly.
In the meantime, this modification procedure can also improve the stability of the scheme by reducing
the magnitude of the first order moments near the discontinuities.
In order to better adapt to our high order scheme, we
choose the HLLC-flux (Harten-Lax-van Leer-contact flux) to be our numerical flux, which is an approximate Riemann
solver by assuming that there are four states in the transition from the left to the right states, thus
staying closer to the real physical situation.

The organization of this paper is as follows:
In Section 2, we describe the reconstruction procedure of the moment-based multi-resolution HWENO scheme for
hyperbolic conservation laws in the one and two dimensions in detail.
In Section 3, we
propose a number of numerical examples to illustrate the accuracy and resolution of our HWENO scheme.
Concluding remarks are given in Section 4.
}

{\section{Moment-based multi-resolution HWENO scheme}}
\label{sec2}
\setcounter{equation}{0}
\setcounter{figure}{0}
\setcounter{table}{0}
{
In this section, we describe the reconstruction procedure of the moment-based multi-resolution HWENO scheme for the one and two dimensional hyperbolic conservation laws, which has sixth order of accuracy in smooth regions and high resolution near
discontinuities.
Here sixth order is simply taken as an example, arbitrarily high order HWENO schemes can be designed following the same lines.
}

{\subsection{One dimensional case}}
{
In this subsection, we first consider the following relatively simple one-dimensional hyperbolic conservation laws
\begin{align}
\label{1D-equ}
\left\{
\begin{aligned}
&u_t+f(u)_x=0,\\
&u(x,0)=u_0(x).
\end{aligned}
\right.
\end{align}

For simplicity, the computational domain is divided by a uniform mesh
$I_i=[x_{i-1/2},x_{i+1/2}]$ with the uniform mesh size $\triangle x=x_{i+1/2}-x_{i-1/2}$,
and the corresponding cell center is denoted by $x_i=\frac{1}{2}(x_{i-1/2}+x_{i+1/2})$.

Firstly, we multiply the governing equation (\ref{1D-equ}) by $\frac{1}{\triangle x}$ and
$\frac{1}{\triangle x}\frac{x-x_i}{\triangle x}$, respectively.
Then, we integrate the resulting equations over the target cell $I_i$ and perform the integration by parts, obtaining the following equations
\begin{align}
\label{1D-equ-int}
\left\{
\begin{aligned}
& \frac{1}{\triangle x}\int_{I_i}u_tdx
=-\frac{1}{\triangle x}\Big[f\big(u(x_{i+1/2},t)\big)
                           -f\big(u(x_{i-1/2},t)\big)\Big],\\
& \frac{1}{\triangle x}\int_{I_i}u_t\frac{x-x_i}{\triangle x}dx
=-\frac{1}{2\triangle x}\Big[f\big(u(x_{i+1/2},t)\big)
                            +f\big(u(x_{i-1/2},t)\big)\Big]
                            +\frac{1}{(\triangle x)^2}\int_{I_i}f(u)dx.
\end{aligned}
\right.
\end{align}
Next, we define the zeroth and first order moments as follows
\begin{align}
\label{1D-moment}
\left\{
\begin{aligned}
&\overline{u}_i(t)=\frac{1}{\triangle x}\int_{I_i}u(x,t)dx,\\
&\overline{v}_i(t)=\frac{1}{\triangle x}\int_{I_i}u(x,t)\frac{x-x_i}{\triangle x}dx,
\end{aligned}
\right.
\end{align}
and then, we swap the spatial integration and time derivation and approximate the flux by a numerical flux, obtaining the following semi-discrete conservative scheme
\begin{align}
\label{1D-equ-int-approximation}
\left\{
\begin{aligned}
&\frac{d\overline{u}_i(t)}{dt}=-\frac{1}{ \triangle x}(\hat{f}_{i+1/2}-\hat{f}_{i-1/2}),\\
&\frac{d\overline{v}_i(t)}{dt}=-\frac{1}{2\triangle x}(\hat{f}_{i+1/2}+\hat{f}_{i-1/2})
                               +\frac{1}{ \triangle x}F_i(u).
\end{aligned}
\right.
\end{align}
Here, the numerical flux is taken to be the HLLC flux
$\hat{f}_{i+1/2}=\hat{f}^{HLLC}(u_{i+1/2}^{-},u_{i+1/2}^{+})$

\begin{align}
\label{f_HLLC}
\hat{f}^{HLLC}(u_L,u_R)=
  \left\{
  \begin{aligned}
   &f_L, \,\,\, 0 \le s_L, \\
   &f_L^{*}, \,\,\, s_L  \le 0 \le s^{*}, \\
   &f_R^{*}, \,\,\, s^{*}\le 0 \le s_R, \\
   &f_R, \,\,\, s_R \le 0,
  \end{aligned}
  \right.
\end{align}
where $f_{L/R}=f(u_{L/R})$, $f_{L/R}^{*}=f_{L/R}+s_{L/R}*(u_{L/R}^{*}-u_{L/R})$; $s_L=\mu_L-c_L*coef_L$, $s_R=\mu_R+c_R*coef_R$,
\begin{align}
\label{COEF}
coef_{L/R}=
  \left\{
  \begin{aligned}
   &1, \,\,\, p^{*} \le p_{L/R}, \\
   &\sqrt{1+\frac{(\gamma+1)*(p^{*}/p_{L/R}-1)}{2\gamma}}, \,\,\, otherwise,
  \end{aligned}
  \right.
\end{align}
\begin{align}
\label{uLR_}
u_{L/R}^{*}=\rho_{L/R}*\frac{s_{L/R}-\mu_{L/R}}{s_{L/R}-s^{*}}*
  \left(
  \begin{array}{c}
   1 \\
   s^{*} \\
   \frac{E_{L/R}}{\rho_{L/R}}+(s^{*}-\mu_{L/R})*\left(s^{*}+\frac{p_{L/R}}{\rho_{L/R}*(s_{L/R}-\mu_{L/R})}\right)
  \end{array}
  \right),
\end{align}
where $p^{*}=\frac{p_L+p_R+(\mu_L-\mu_R)*(c^{*}*\rho^{*})}{2}$, $s^{*}=\frac{\mu_L+\mu_R+(p_L-p_R)/(c^{*}*\rho^{*})}{2}$; $\rho^{*}=\frac{\rho_L+\rho_R}{2}$, $c^{*}=\frac{c_L+c_R}{2}$.
The integral term $F_i(u)$ is approximated by a four-point Gauss-Lobatto integration
\begin{equation}
\label{F}
F_i(u)=\frac{1}{\triangle x}\int_{I_i}f(u)dx\approx\sum_{l=1}^4 \omega_l f\big(u(x_l^{GL},t)\big),
\end{equation}
where the weights are $\omega_1=\omega_4=\frac{1}{12}$, $\omega_2=\omega_3=\frac{5}{12}$ and the quadrature points on the target cell $I_i$ are
\begin{gather*}
x_1^{GL}=x_{i-1/2},x_2^{GL}=x_{i-\sqrt{5}/10},x_3^{GL}=x_{i+\sqrt{5}/10},x_4^{GL}=x_{i+1/2},
\end{gather*}
where $x_{i+\xi}=x_{i}+\xi \triangle x$.

Finally, we present the spatial reconstruction procedure of the point values $\{u_{i-1/2}^{+}$, $u_{i-\sqrt{5}/10}$, $u_{i+\sqrt{5}/10}$, $u_{i+1/2}^{-}\}$ from the given moment values $\{\overline{u},\overline{v}|_{i-1};\overline{u},\overline{v}|_{i};\overline{u},\overline{v}|_{i+1}\}$ in detail as follows:

{\bf{\em The 1D Reconstruction Algorithm:}}

{\bf{\em Step 1.}} Reconstruct the Gauss-Lobatto point values of $u$.

{\bf{\em Step 1.1.}} Reconstruct a series of polynomials of different degrees.

First, we select a series of central spatial stencils and reconstruct a zeroth degree polynomial $q_1(x)$, a quadratic polynomial $q_2(x)$, a cubic polynomial $q_3(x)$ and a quintic polynomial $q_4(x)$, respectively, such that
\begin{align}
\label{q}
\left.
\begin{aligned}
      \frac{1}{\triangle x}\int_{I_ k   }q_1(x)                             dx
    &=\overline{u}_ k   ,   k   =    i    ; \\
      \frac{1}{\triangle x}\int_{I_ k   }q_2(x)                             dx
    &=\overline{u}_ k   ,   k   =i-1,i,i+1; \\
      \frac{1}{\triangle x}\int_{I_ k   }q_3(x)                             dx
    &=\overline{u}_ k   ,   k   =i-1,i,i+1;
\quad \frac{1}{\triangle x}\int_{I_{k_x}}q_3(x)\frac{x-x_{k_x}}{\triangle x}dx
     =\overline{v}_{k_x},  {k_x}=    i    ; \\
      \frac{1}{\triangle x}\int_{I_ k   }q_4(x)                             dx
    &=\overline{u}_ k   ,   k   =i-1,i,i+1;
\quad \frac{1}{\triangle x}\int_{I_{k_x}}q_4(x)\frac{x-x_{k_x}}{\triangle x}dx
     =\overline{v}_{k_x},  {k_x}=i-1,i,i+1.
\end{aligned}
\right.
\end{align}

Then, we obtain the following equivalent expressions for the above reconstructed polynomials
\begin{align}
\label{p}
p_{l_2}(x)=
  \left\{
  \begin{aligned}
   &q_1(x), \,\,\, l_2=1, \\
   &\frac{1}{\gamma_{l_2,l_2}}q_{l_2}(x)-\sum_{l=1}^{l_2-1}\frac{\gamma_{l,l_2}}{\gamma_{l_2,l_2}}p_l(x),
   \,\,\, l_2=2,3,4,
  \end{aligned}
  \right.
\end{align}
with $\sum_{l=1}^{l_2}\gamma_{l,l_2}=1,\gamma_{l_2,l_2}\neq0,l_2=2,3,4$, where these
$\gamma_{l_1,l_2}$ for $l_1=1,...,l_2;l_2=2,3,4$ are the linear weights and are defined as
\begin{equation}
\label{linear-weight}
\gamma_{l_1,l_2}=\frac{\overline{\gamma}_{l_1,l_2}}{\sum_{l=1}^{l_2}\overline{\gamma}_{l,l_2}};
\quad \overline{\gamma}_{l_1,l_2}=10^{l_1-1}; \quad l_1=1,...,l_2; \,\, l_2=2,3,4.
\end{equation}
Putting these linear weights into (\ref{p}), we obtain the following relations
\begin{align}
\label{pq}
\left.
\begin{aligned}
p_1(x)&=                 q_1(x)                       , \\
p_2(x)&=\frac{  11}{10  }q_2(x)-\frac{  1}{10  }q_1(x), \\
p_3(x)&=\frac{ 111}{100 }q_3(x)-\frac{ 11}{100 }q_2(x), \\
p_4(x)&=\frac{1111}{1000}q_4(x)-\frac{111}{1000}q_3(x).
\end{aligned}
\right.
\end{align}

{\bf{\em Step 1.2.}} Compute the corresponding nonlinear weights of the above polynomials.

First, we compute the smoothness indicator $\beta_{l_2}$ of function $p_{l_2}(x)$ in the interval $I_i$:
\begin{equation}
\label{beta}
\beta_{l_2}=\sum_{\alpha=1}^{\kappa}\int_{I_i}{\triangle x}^{2\alpha-1}
\left(\frac{d^{\alpha}p_{l_2}(x)}{dx^{\alpha}}\right)^2dx, \quad l_2=2,3,4,
\end{equation}
where $\kappa=2,3,5$ for $l_2=2,3,4$.
It is important to note that the definition of $\beta_{1}$ is different, where a new  polynomial $q_1^{new}(x)$ is required and is defined as follows:
\newline
(1) We reconstruct two polynomials $q_{1L}(x)$ and $q_{1R}(x)$, such that
\begin{align}
\label{q1-LR}
\left.
\begin{aligned}
      \frac{1}{\triangle x}\int_{I_k}q_{1L}(x)dx
    &=\overline{u}_k, \,\,\, k=i-1,i; \\
      \frac{1}{\triangle x}\int_{I_k}q_{1R}(x)dx
    &=\overline{u}_k, \,\,\, k=i,i+1,
\end{aligned}
\right.
\end{align}
and then, we obtain their associated smoothness indicators
\begin{equation}
\label{beta1-LR}
      \beta_{1L}=(\overline{u}_{i}-\overline{u}_{i-1})^2,
\quad \beta_{1R}=(\overline{u}_{i+1}-\overline{u}_{i})^2,
\end{equation}
and the absolute difference $\tau_{1}$ between $\beta_{1L}$ and $\beta_{1R}$
\begin{equation}
\label{tau1}
\tau_{1}=|\beta_{1R}-\beta_{1L}|^{2},
\end{equation}
where the selection of the power is to be consistent with the definition (\ref{tau4}) of $\tau_{4}$ later.
\newline
(2) We give these two polynomials $q_{1L}(x)$ and $q_{1R}(x)$ the same linear weights $\gamma_{1L}=\gamma_{1R}=\frac{1}{2}$ and calculate the corresponding nonlinear weights as
\begin{equation}
\label{w1-LR}
      \omega_{1L}=\frac{\overline{\omega}_{1L}}
{\overline{\omega}_{1L}+\overline{\omega}_{1R}},
\quad \omega_{1R}=\frac{\overline{\omega}_{1R}}
{\overline{\omega}_{1L}+\overline{\omega}_{1R}},
\end{equation}
\begin{equation}
\label{w1_-LR}
      \overline{\omega}_{1L}=\gamma_{1L}
\left(1+\frac{\tau_{1}}{\beta_{1L}+\varepsilon}\right),
\quad \overline{\omega}_{1R}=\gamma_{1R}
\left(1+\frac{\tau_{1}}{\beta_{1R}+\varepsilon}\right),
\end{equation}
where $\varepsilon=10^{-10}$ is applied to avoid the denominator of (\ref{w1_-LR}) to be zero.
\newline
(3) We obtain a new polynomial
\begin{equation}
\label{q1_new}
q_1^{new}(x)=\omega_{1L}q_{1L}(x)+\omega_{1R}q_{1R}(x),
\end{equation}
and set $\beta_1$ to be
\begin{equation}
\label{beta1}
\beta_1=\int_{I_i}{\triangle x}\left(\frac{d q_1^{new}(x)}{dx}\right)^2dx
       =\Big(\omega_{1L}\big(\overline{u}_{i}-\overline{u}_{i-1}\big)
       +     \omega_{1R}\big(\overline{u}_{i+1}-\overline{u}_{i}\big)\Big)^2.
\end{equation}

Then, we still adopt the idea of WENO-Z as shown in \cite{ccd2011}
to define the corresponding nonlinear weights
\begin{equation}
\label{non-linear-weight}
\omega_{l_1,4}=\frac{\overline{\omega}_{l_1,4}}
{\sum_{l=1}^{4}\overline{\omega}_{l,4}},
\quad \overline{\omega}_{l_1,4}=\gamma_{l_1,4}
\left(1+\frac{\tau_{4}}{\beta_{l_1}+\varepsilon}\right), \quad l_1=1,...,4,
\end{equation}
where $\varepsilon$ is also taken to be $10^{-10}$ and the quantity $\tau_{4}$ is defined to be the absolute difference among above smoothness indicators
\begin{equation}
\label{tau4}
\tau_{4}=\left(\frac{\sum_{l=1}^{3}|\beta_{4}-\beta_{l}|}{3}\right)^{2}.
\end{equation}

{\bf{\em Step 1.3.}} Obtain an approximation polynomial $u_i(x)$ of $u(x)$.

The new reconstruction polynomial $u_i(x)$ of $u(x)$ is defined as
\begin{equation}
\label{re-p}
u_i(x)=\sum_{l=1}^{4}\omega_{l,4}p_l(x),
\end{equation}
and the Gauss-Lobatto point values that we need are taken to be
\begin{equation}
\label{re-u}
u_{i-1/2}^{+}=u_i(x_{i-1/2}),
u_{i\mp\sqrt{5}/10}=u_i(x_{i\mp\sqrt{5}/10}),
u_{i+1/2}^{-}=u_i(x_{i+1/2}).
\end{equation}

{\bf{\em Step 2.}} Update the Gauss-Lobatto point values of $u$ in the troubled-cells.

{\bf{\em Step 2.1.}} Identify the troubled-cells.

The so-called troubled-cells are those cells that may contain discontinuities. In 2005, Qiu and Shu systematically
investigated and compared a few troubled-cell indicators for the Runge-Kutta discontinuous Galerkin method
in \cite{qs2005}. Here, we choose the KXRCF troubled-cell indicator proposed by Krivodonova et al. in \cite{kxrcf2004}
to identify the troubled-cells, and its judgment criterion is that the target cell $I_i$ is identified to be a
troubled-cell, if
\begin{equation}
\label{KXRCF}
\aleph_{i}=\frac{\big|\int_{\partial I_i^{-}}\big(u_i(x)-u_{n_i}(x)\big)ds\big|}
{h_i^{\frac{l+1}{3}}\big|\partial I_i^{-}\big| ||u_i(x)||}>1,
\end{equation}
where $\partial I_i^{-}$ is the inflow boundary ($\overrightarrow{v}\cdot \overrightarrow{n}<0$, $\overrightarrow{v}$ is the velocity of the flow and $\overrightarrow{n}$ is the outer normal vector to $\partial I_i$), $I_{n_i}$ is the neighbor of $I_{i}$ on the side of $\partial I_i^{-}$, $h_i$ is the length of the cell $I_{i}$, the parameter $l$ (i.e. the degree of $u_i(x)$) is taken to be 5, $u_i(x)$ is the approximation polynomial of $u(x)$ obtained in {\bf{\em Step 1}} above and the norm is
taken to be the $L^{\infty}$ norm.

{\bf{\em Step 2.2.}} Modify the first order moments in the troubled-cells.

If the target cell $I_i$ is identified to be a troubled-cell, we would like to modify the first order moment $\overline{v}_i$ in it. First, we reconstruct a quartic polynomial, which satisfies
\begin{equation}
\label{p0}
       \frac{1}{\triangle x}\int_{I_{k  }}p_0(x)                             dx
       =\overline{u}_{k  }, \quad k  =i-1,i,i+1;
\quad  \frac{1}{\triangle x}\int_{I_{k_x}}p_0(x)\frac{x-x_{k_x}}{\triangle x}dx
       =\overline{v}_{k_x}, \quad k_x=i-1,  i+1,
\end{equation}
and then, we modify the first order moment $\overline{v}_i$ as
\begin{equation}
\label{vi}
\overline{v}_i=
\frac{1}{\triangle x}\int_{I_i}p_0(x)\frac{x-x_i}{\triangle x}dx=
-\frac{ 5}{76}\overline{u}_{i-1}+\frac{ 5}{76}\overline{u}_{i+1}
-\frac{11}{38}\overline{v}_{i-1}-\frac{11}{38}\overline{v}_{i+1}.
\end{equation}

{\bf{\em Step 2.3.}} Update the Gauss-Lobatto point values of $u$ for these troubled-cells.

After modifying the first order moments in the troubled-cells, we repeat the reconstruction process {\bf{\em Step 1}}
for these troubled-cells to update corresponding Gauss-Lobatto point values of $u$.

{\bf{\em Step 3.}} Discretize the semi-discrete scheme in time.

After all these Gauss-Lobatto point values are obtained, we substitute them into the formula of the numerical flux.
Then, we discretize (\ref{1D-equ-int-approximation}) by a third-order TVD Runge-Kutta method in time
\begin{align}
\label{Runge-Kutta}
\left\{
\begin{aligned}
&U^{(1)}=U^{n}+\Delta t L(U^{n}),\\
&U^{(2)}=\frac{3}{4}U^{n}+\frac{1}{4}U^{(1)}+\frac{1}{4}\Delta t L(U^{(1)}),\\
&U^{n+1}=\frac{1}{3}U^{n}+\frac{2}{3}U^{(2)}+\frac{2}{3}\Delta t L(U^{(2)}),
\end{aligned}
\right.
\end{align}
to obtain a fully discrete scheme.

{\bf{\em Remark 1.}}
In {\bf{\em Step 1}} above, through a series of Taylor expansion analyses, we can verify that $\beta_{l}=u^{'2}\triangle x^2+O(\triangle x^4), l=1,2,3,4$, thus $\beta_{4}-\beta_{l}=O(\triangle x^4)$ for $l<4$ and $\tau_{4}=O(\triangle x^8)$, then \begin{equation}
\label{error-analysis}
\begin{split}
u(x)-u_i(x)&=u(x)-\sum_{l=1}^{4}\omega_{l,4}p_l(x) \\
           &=\left[u(x)-\sum_{l=1}^{4}\gamma_{l,4}p_l(x)\right]
            +\sum_{l=1}^{4}(\omega_{l,4}-\gamma_{l,4})\Big(u\big(x\big)-p_l\big(x\big)\Big) \\
           &=O(\triangle x^{6})+O(\triangle x^{6})*O(\triangle x) \\
           &=O(\triangle x^{6}).
\end{split}
\end{equation}
According to above Taylor expansion analyses, we can prove that our reconstruction algorithm can obtain sixth order of accuracy in
smooth regions.

When there is a discontinuity in the target cell $I_i$, its first order moment will become $\overline{v}_i=O(1)$, thus $\beta_{1}=\beta_{2}=O(\triangle x^2), \beta_{3}=\beta_{4}=O(1)$ and  $\tau_{4}=O(1)$, then
\begin{align}
\label{dis-error-analysis}
\left.
\begin{aligned}
&\overline{\omega}_{1,4}=\gamma_{1,4}\big(1+O(\triangle x^{-2})\big), \quad
 \overline{\omega}_{2,4}=\gamma_{2,4}\big(1+O(\triangle x^{-2})\big);\\
&\overline{\omega}_{3,4}=\gamma_{3,4}\big(1+O(1)\big), \quad
 \overline{\omega}_{4,4}=\gamma_{4,4}\big(1+O(1)\big).
\end{aligned}
\right.
\end{align}
That is to say, the proportion of the first two layers will be magnified slightly, which will cause the discontinuity to be smoothed out to a certain extent.
But if we modify the first order moment of the troubled-cell $I_i$, its first order moment will become $\overline{v}_i=O(\triangle x)$, thus $\beta_{l}=u^{'2}\triangle x^2+O(\triangle x^4), l=1,2,3,4$, and  $\tau_{4}=O(\triangle x^8)$, then $\overline{\omega}_{l,4}=\gamma_{l,4}\big(1+O(\triangle x^{6})\big), l=1,2,3,4$, i.e. the nonlinear weight of each layer is close to their corresponding linear weight, thus the resolution near the discontinuities can be increased significantly.

As in our previous paper, the choice of these linear weights is also not unique.
Even though different choices of the linear weights do not affect the order of accuracy in smooth regions, they do
affect the resolution near discontinuities.
That is to say, the bigger the linear weights are for higher-degree polynomials, the steeper the shock transitions become,
but also the more unstable the scheme becomes. This modification procedure can also improve the stability of the scheme by reducing the magnitude of the first order moments near the discontinuities.

{\bf{\em Remark 2.}}
In {\bf{\em Step 2}} above, we can also choose other indicators to identify troubled-cells, such as
the minmod-based total variation bounded (TVB) limiter in \cite{cs1989},
                      moment limiter of Biswas, Devine and Flaherty in \cite{bdf1994},
a modification of the moment limiter by Burbeau, Sagaut and Bruneau in \cite{bsb2001},
the monotonicity-preserving(MP) limiter in \cite{sh1997},
a modification of the MP limiter in \cite{rm2001},
a troubled-cell indicator based on Harten's subcell resolution idea in \cite{h1989}
and so on.
But as shown in Section 3, the KXRCF troubled-cell indicator works pretty well for our scheme in the one-dimensional case.
What needs a special attention is that,
for the one dimensional scalar equation,
     the solution $u$ is defined as  our indicator variable, and then the corresponding
     $\overrightarrow{v}=f^{'}(u)$;
for the one dimensional Euler system,
only the density $\rho$ is set to be our indicator variable, and then the corresponding
     $\overrightarrow{v}=\mu$ is the velocity of the fluid.
In short, for the one dimensional case, the line integral average in the formula (\ref{KXRCF}) is actually the boundary point value, that is
\begin{equation}
\label{line-integral-average}
\frac{1}{\big|\partial I_i^{-}\big|}\left|\int_{\partial I_i^{-}}\big(u_i(x)-u_{n_i}(x)\big)ds\right|
=\left|
  (u_{i-\frac{1}{2}}^{+}-u_{i-\frac{1}{2}}^{-})*sf( \overrightarrow{v}_{i-\frac{1}{2}})
 +(u_{i+\frac{1}{2}}^{-}-u_{i+\frac{1}{2}}^{+})*sf(-\overrightarrow{v}_{i+\frac{1}{2}})
 \right|,
\end{equation}

and the norm $||u_i(x)||$ is
taken to be the maximum norm of all the Gauss-Lobatto point values in the cell $I_i$
(i.e. $||u_i(x)|| \approx \max \{ |u_{i-1/2}^{+}|, |u_{i-\sqrt{5}/10}|, |u_{i+\sqrt{5}/10}|, |u_{i+1/2}^{-}|\}$).
Note that all the values used in the troubled-cell indicator are already obtained in the reconstruction process {\bf{\em Step 1}},
thus there is no need to reconstruct an extra polynomial as in \cite{zcq2020,zq2020}.

According to \cite{kxrcf2004},
\begin{align}
\label{analysis1}
  \int_{\partial I_i^{-}}\big(u_i-u_{n_i}\big)ds
&=\int_{\partial I_i^{-}}\big(u_i-u\big)ds+\int_{\partial I_{n_i}^{+}}\big(u-u_{n_i}\big)ds\\
&=\notag
  \left\{
  \begin{aligned}
   &O(h^{p+2})+O(h^{2p+2})=O(h^{p+2}), u|_{\partial I_i} \ is \ smooth, \\
   &O(h)+O(h)=O(h), u|_{\partial I_i} \ is \ discontinuous.
  \end{aligned}
  \right.
\end{align}
Here, we let the troubled-cell indicator converge to the smooth case twice as fast as the discontinuous case, then as either $h\rightarrow 0$ or $p\rightarrow \infty$
\begin{align}
\label{analysis2}
  \aleph_{i}=
  \left\{
  \begin{aligned}
   &O(h^{\frac{2(l+1)}{3}})       \stackrel{}{\longrightarrow}0     ,
   u|_{\partial I_i} \ is \ smooth, \\
   &\frac{1}{O(h^{\frac{l+1}{3}})}\stackrel{}{\longrightarrow}\infty,
   u|_{\partial I_i} \ is \ discontinuous,
  \end{aligned}
  \right.
\end{align}
thus the troubled-cell indicator can be defined as

\begin{align}
\label{analysis3}
  \left\{
  \begin{aligned}
   &\aleph_{i}<1, \,\,\, I_i \mbox{ is not a troubled-cell}, \\
   &\aleph_{i}>1, \,\,\, I_i \mbox{ is a troubled-cell}.
  \end{aligned}
  \right.
\end{align}
}

{\subsection{Two dimensional case}}
{
In this subsection, we then consider the following more complicated two-dimensional hyperbolic conservation laws
\begin{align}
\label{2D-equ}
\left\{
\begin{aligned}
&u_t+f(u)_x+g(u)_y=0,\\
&u(x,y,0)=u_0(x,y).
\end{aligned}
\right.
\end{align}

Still for the sake of simplicity, the computational domain is divided by a uniform mesh $I_{i,j}=$ $[x_{i-1/2},x_{i+1/2}]\times[y_{j-1/2},y_{j+1/2}]$ with the uniform mesh sizes
$\triangle x=x_{i+1/2}-x_{i-1/2}$ in the $x$-direction and $\triangle y=y_{j+1/2}-y_{j-1/2}$ in the $y$-direction, and the corresponding cell center is denoted by $(x_i,y_j)
=\Big(\frac{1}{2}\big(x_{i-1/2}+x_{i+1/2}\big),\frac{1}{2}\big(y_{j-1/2}+y_{j+1/2}\big)\Big)$.

Firstly, we multiply the governing equation (\ref{2D-equ}) by $\frac{1}{\triangle x\triangle y}$,
$\frac{1}{\triangle x\triangle y}\frac{x-x_i}{\triangle x}$ and $\frac{1}{\triangle x\triangle y}\frac{y-y_j}{\triangle y}$,
respectively.
Then, we integrate the resulting equations over the target cell $I_{i,j}$ and perform the integration by parts, obtaining the following equations
\begin{align}
\label{2D-equ-int}
\left\{
\begin{aligned}
& \frac{1}{\triangle x\triangle y}\int\limits_{I_{i,j}}u_tdxdy
=-\frac{1}{\triangle x\triangle y}\int\limits_{y_{j-1/2}}^{y_{j+1/2}}
  \Big[f\big(u(x_{i+1/2},y,t)\big)-f\big(u(x_{i-1/2},y,t)\big)\Big]dy\\
&-\frac{1}{\triangle x\triangle y}\int\limits_{x_{i-1/2}}^{x_{i+1/2}}
  \Big[g\big(u(x,y_{j+1/2},t)\big)-g\big(u(x,y_{j-1/2},t)\big)\Big]dx,\\
& \frac{1}{ \triangle x   \triangle y}\int\limits_{I_{i,j}}u_t\frac{x-x_i}{\triangle x}dxdy
=-\frac{1}{2\triangle x   \triangle y}\int\limits_{y_{j-1/2}}^{y_{j+1/2}}
  \Big[f\big(u(x_{i+1/2},y,t)\big)+f\big(u(x_{i-1/2},y,t)\big)\Big]dy\\
&+\frac{1}{(\triangle x)^2\triangle y}\int\limits_{I_{i,j}}f(u)dxdy
 -\frac{1}{ \triangle x   \triangle y}\int\limits_{x_{i-1/2}}^{x_{i+1/2}}
  \Big[g\big(u(x,y_{j+1/2},t)\big)-g\big(u(x,y_{j-1/2},t)\big)\Big]
  \frac{x-x_i}{\triangle x}dx,\\
& \frac{1}{ \triangle x   \triangle y}\int\limits_{I_{i,j}}u_t\frac{y-y_j}{\triangle y}dxdy
=-\frac{1}{ \triangle x   \triangle y}\int\limits_{y_{j-1/2}}^{y_{j+1/2}}
  \Big[f\big(u(x_{i+1/2},y,t)\big)-f\big(u(x_{i-1/2},y,t)\big)\Big]
  \frac{y-y_j}{\triangle y}dy\\
&-\frac{1}{2\triangle x   \triangle y}\int\limits_{x_{i-1/2}}^{x_{i+1/2}}
  \Big[g\big(u(x,y_{j+1/2},t)\big)+g\big(u(x,y_{j-1/2},t)\big)\Big]dx
 +\frac{1}{\triangle x(\triangle y)^2}\int\limits_{I_{i,j}}g(u)dxdy.
\end{aligned}
\right.
\end{align}
Next, we define the zeroth and first order moments in the $x$ and $y$ directions as follows
\begin{align}
\label{2D-moment}
\left\{
\begin{aligned}
&\tilde{\overline{u}}_{i,j}(t)=\frac{1}{\triangle x\triangle y}\int_{I_{i,j}}u(x,y,t)dxdy,\\
&\tilde{\overline{v}}_{i,j}(t)=\frac{1}{\triangle x\triangle y}
\int_{I_{i,j}}u(x,y,t)\frac{x-x_i}{\triangle x}dxdy,\\
&\tilde{\overline{w}}_{i,j}(t)=\frac{1}{\triangle x\triangle y}
\int_{I_{i,j}}u(x,y,t)\frac{y-y_j}{\triangle y}dxdy,
\end{aligned}
\right.
\end{align}
and then, we swap the spatial integration and time derivation and approximate the flux by
a numerical flux, obtaining the following semi-discrete conservative scheme

\begin{align}
\label{2D-equ-int-approximation}
\left\{
\begin{aligned}
& \frac{d\tilde{\overline{u}}_{i,j}(t)}{dt}
=-\frac{1}{\triangle x\triangle y}\int_{y_{j-1/2}}^{y_{j+1/2}}
  \Big[\hat{f}\big(u(x_{i+1/2},y)\big)-\hat{f}\big(u(x_{i-1/2},y)\big)\Big]dy\\
&-\frac{1}{\triangle x\triangle y}\int_{x_{i-1/2}}^{x_{i+1/2}}
  \Big[\hat{g}\big(u(x,y_{j+1/2})\big)-\hat{g}\big(u(x,y_{j-1/2})\big)\Big]dx,\\
& \frac{d\tilde{\overline{v}}_{i,j}(t)}{dt}
=-\frac{1}{2\triangle x\triangle y}\int_{y_{j-1/2}}^{y_{j+1/2}}
  \Big[\hat{f}\big(u(x_{i+1/2},y)\big)+\hat{f}\big(u(x_{i-1/2},y)\big)\Big]dy\\
&+\frac{1}{ \triangle x}F_{i,j}(u)
 -\frac{1}{ \triangle x\triangle y}\int_{x_{i-1/2}}^{x_{i+1/2}}
  \Big[\hat{g}\big(u(x,y_{j+1/2})\big)-\hat{g}\big(u(x,y_{j-1/2})\big)\Big]
  \frac{x-x_i}{\triangle x}dx,\\
& \frac{d\tilde{\overline{w}}_{i,j}(t)}{dt}
=-\frac{1}{ \triangle x\triangle y}\int_{y_{j-1/2}}^{y_{j+1/2}}
  \Big[\hat{f}\big(u(x_{i+1/2},y)\big)-\hat{f}\big(u(x_{i-1/2},y)\big)\Big]
  \frac{y-y_j}{\triangle y}dy\\
&-\frac{1}{2\triangle x\triangle y}\int_{x_{i-1/2}}^{x_{i+1/2}}
  \Big[\hat{g}\big(u(x,y_{j+1/2})\big)+\hat{g}\big(u(x,y_{j-1/2})\big)\Big]dx
 +\frac{1}{ \triangle y}G_{i,j}(u).
\end{aligned}
\right.
\end{align}
Here the numerical fluxes are still taken to be the HLLC fluxes,
and the integral terms are also approximated by a four-point Gauss-Lobatto integration, for instance
\begin{equation}
\label{integral1}
F_{i,j}(u)=       \frac{1}{\triangle x\triangle y}\int_{I_{i,j}}f(u)dxdy
          \approx \sum_{k=1}^4\sum_{l=1}^4 \omega_k\omega_l f\big(u(x_{k}^{GL},y_{l}^{GL})\big),
\end{equation}
\begin{equation}
\label{integral2}
\frac{1}{\triangle y}\int_{y_{j-1/2}}^{y_{j+1/2}}\hat{f}\big(u(x_{i+1/2},y)\big)dy
\approx \sum_{l=1}^4 \omega_l\hat{f}\big(u(x_{i+1/2},y_{l}^{GL})\big),
\end{equation}
\begin{equation}
\label{integral3}
\frac{1}{\triangle x}\int_{x_{i-1/2}}^{x_{i+1/2}}
\frac{x-x_i}{\triangle x}\hat{g}\big(u(x,y_{j+1/2})\big)dx
\approx \sum_{k=1}^4 \omega_k
\frac{x_{k}^{GL}-x_i}{\triangle x}\hat{g}\big(u(x_{k}^{GL},y_{j+1/2})\big),
\end{equation}
where the weights are $\omega_1=\omega_4=\frac{1}{12}$, $\omega_2=\omega_3=\frac{5}{12}$ and the quadrature points on the target cell $I_{i,j}$ are
\begin{gather*}
x_1^{GL}=x_{i-1/2},x_2^{GL}=x_{i-\sqrt{5}/10},x_3^{GL}=x_{i+\sqrt{5}/10},x_4^{GL}=x_{i+1/2};\\
y_1^{GL}=y_{j-1/2},y_2^{GL}=y_{j-\sqrt{5}/10},y_3^{GL}=y_{j+\sqrt{5}/10},y_4^{GL}=y_{j+1/2},
\end{gather*}
where $x_{i+\xi}=x_{i}+\xi \triangle x$ and $y_{j+\eta}=y_{j}+\eta \triangle y$.

Finally, we present the spatial reconstruction procedure of the point values $\{u_{i\mp1/2,j+\eta_{l}}^{\pm}|l=1,2,3,4\}$,
$\{u_{i+\xi_{k},j\mp1/2}^{\pm}|k=1,2,3,4\}$ and
$\{u_{i+\xi_{k},j+\eta_{l}}|k=2,3;l=2,3\}$
($\xi_{1}=\eta_{1}=-1/2$, $\xi_{2}=\eta_{2}=-\sqrt{5}/10$, $\xi_{3}=\eta_{3}=\sqrt{5}/10$ and $\xi_{4}=\eta_{4}=1/2$) from the given cell-average values $\{\tilde{\overline{u}},\tilde{\overline{v}},\tilde{\overline{w}}|_{k_x,k_y};k_x=i-1,i,i+1;k_y=j-1,j,j+1\}$ in detail as follows:

{\bf{\em The 2D Reconstruction Algorithm:}}

{\bf{\em Step 1.}} Reconstruct the Gauss-Lobatto point values of $u$.

\begin{center}
\begin{tabular}{cccc}
\cline{1-3}
\multicolumn{1}{|c|}{7} & \multicolumn{1}{c|}{8} & \multicolumn{1}{c|}{9} & $j+1$ \\ \cline{1-3}
\multicolumn{1}{|c|}{4} & \multicolumn{1}{c|}{5} & \multicolumn{1}{c|}{6} & $j  $ \\ \cline{1-3}
\multicolumn{1}{|c|}{1} & \multicolumn{1}{c|}{2} & \multicolumn{1}{c|}{3} & $j-1$ \\ \cline{1-3}
 $i-1$                  &  $i$                   &  $i+1$                 &
\end{tabular}
\\ The big stencil and its new labels.
\end{center}

{\bf{\em Step 1.1.}} Reconstruct a series of polynomials of different degrees.

{\bf{\em Step 1.1.1.}} Reconstruct a zeroth degree polynomial $q_1(x,y)$ such that
\begin{equation}
\label{2D-q1}
\frac{1}{\triangle x\triangle y}\int_{I_k}q_1(x,y)dxdy=\tilde{\overline{u}}_{k},
\quad k=5.
\end{equation}

{\bf{\em Step 1.1.2.}} Reconstruct a quadratic polynomial $q_2(x,y)$ such that
\begin{gather}
\label{2D-q2}
\frac{1}{\triangle x\triangle y}\int_{I_k}q_2(x,y)dxdy=\tilde{\overline{u}}_{k},
\quad k=1,...,9.
\end{gather}

{\bf{\em Step 1.1.3.}} Reconstruct a cubic polynomial $q_3(x,y)$ such that
\begin{gather}
\label{2D-q3}
\notag \frac{1}{\triangle x\triangle y}\int_{I_ k   }q_3(x,y)                             dxdy
=\tilde{\overline{u}}_{k  },\quad k   =1,...,9;\\
       \frac{1}{\triangle x\triangle y}\int_{I_{k_x}}q_3(x,y)\frac{x-x_{k_x}}{\triangle x}dxdy
=\tilde{\overline{v}}_{k_x},\quad{k_x}=5;      \\
\notag \frac{1}{\triangle x\triangle y}\int_{I_{k_y}}q_3(x,y)\frac{y-y_{k_y}}{\triangle y}dxdy
=\tilde{\overline{w}}_{k_y},\quad{k_y}=5.
\end{gather}

{\bf{\em Step 1.1.4.}} Reconstruct a quintic polynomial $q_4(x,y)$ such that
\begin{gather}
\label{2D-q4}
\notag \frac{1}{\triangle x\triangle y}\int_{I_ k   }q_4(x,y)                             dxdy
=\tilde{\overline{u}}_{k  },\quad k   =1,...,9;      \\
       \frac{1}{\triangle x\triangle y}\int_{I_{k_x}}q_4(x,y)\frac{x-x_{k_x}}{\triangle x}dxdy
=\tilde{\overline{v}}_{k_x},\quad{k_x}=1,3,4,5,6,7,9;\\
\notag \frac{1}{\triangle x\triangle y}\int_{I_{k_y}}q_4(x,y)\frac{y-y_{k_y}}{\triangle y}dxdy
=\tilde{\overline{w}}_{k_y},\quad{k_y}=1,2,3,5,7,8,9.
\end{gather}
What needs special attention is that
the number of equations is greater than the number of unknowns when we reconstruct the quadratic polynomial $q_2(x,y)$, the cubic polynomial $q_3(x,y)$ and
the quintic polynomial $q_4(x,y)$.
To solve this problem, we require these polynomials must have the same cell average as $u$ on the target cell $I_{i,j}$
(to maintain conservation) and match the other conditions in a least square sense as described in \cite{hs1999}.

{\bf{\em Step 1.1.5.}} Then, we further manipulate the above reconstructed polynomials to obtain the following equivalent expressions
\begin{align}
\label{2D-p}
p_{l_2}(x,y)=
  \left\{
  \begin{aligned}
   &q_1(x,y), \,\,\, l_2=1, \\
   &\frac{1}{\gamma_{l_2,l_2}}q_{l_2}(x,y)-\sum_{l=1}^{l_2-1}\frac{\gamma_{l,l_2}}{\gamma_{l_2,l_2}}p_l(x,y),
   \,\,\, l_2=2,3,4,
  \end{aligned}
  \right.
\end{align}
with $\sum_{l=1}^{l_2}\gamma_{l,l_2}=1,\gamma_{l_2,l_2}\neq0,l_2=2,3,4$,
where these $\gamma_{l_1,l_2}$ for $l_1=1,...,l_2;l_2=2,3,4$ are still the linear weights and are defined as (\ref{linear-weight}).
Likewise putting these linear weights into (\ref{2D-p}), we obtain the following relations
\begin{align}
\label{2D-pq}
\left.
\begin{aligned}
p_1(x,y)&=                 q_1(x,y)                         , \\
p_2(x,y)&=\frac{  11}{10  }q_2(x,y)-\frac{  1}{10  }q_1(x,y), \\
p_3(x,y)&=\frac{ 111}{100 }q_3(x,y)-\frac{ 11}{100 }q_2(x,y), \\
p_4(x,y)&=\frac{1111}{1000}q_4(x,y)-\frac{111}{1000}q_3(x,y).
\end{aligned}
\right.
\end{align}

{\bf{\em Step 1.2.}} Compute the corresponding nonlinear weights of the above polynomials.

First, we compute the smoothness indicator $\beta_{l_2}$ of the function $p_{l_2}(x,y)$ in the interval $I_{i,j}$:
\begin{equation}
\label{2D-beta}
\beta_{l_2}=\sum_{|\alpha|=1}^{\kappa}\int_{I_{i,j}}|I_{i,j}|^{|\alpha|-1}
\left(\frac{\partial^{|\alpha|}}
{\partial x^{{\alpha}_x}\partial y^{{\alpha}_y}}p_{l_2}(x,y)\right)^2dxdy, \quad l_2=2,3,4,
\end{equation}
where $\alpha=(\alpha_x,\alpha_y)$, $|\alpha|=\alpha_x+\alpha_y$ and $\kappa=2,3,5$ for $l_2=2,3,4$. Note that the definition of $\beta_{1}$ is an exception, where a new polynomial $q_1^{new}(x,y)$ is required and is defined as follows:
\newline
(1) We reconstruct four polynomials $q_{1k}(x,y)$ for $k=1,2,3,4$, such that
\begin{align}
\label{2D-q1-1234}
\left.
\begin{aligned}
  \frac{1}{\triangle x\triangle y}\int_{I_k}q_{11}(x,y)dxdy
&=\tilde{\overline{u}}_{k}, \quad k=4,5,8;\\
  \frac{1}{\triangle x\triangle y}\int_{I_k}q_{12}(x,y)dxdy
&=\tilde{\overline{u}}_{k}, \quad k=5,6,8;\\
  \frac{1}{\triangle x\triangle y}\int_{I_k}q_{13}(x,y)dxdy
&=\tilde{\overline{u}}_{k}, \quad k=2,5,6;\\
  \frac{1}{\triangle x\triangle y}\int_{I_k}q_{14}(x,y)dxdy
&=\tilde{\overline{u}}_{k}, \quad k=2,4,5,
\end{aligned}
\right.
\end{align}
and then, we obtain their associated smoothness indicators
\begin{align}
\label{2D-beta1-1234}
\left.
\begin{aligned}
\beta_{11}&=(\tilde{\overline{u}}_{i,j}-\tilde{\overline{u}}_{i-1,j})^2
           +(\tilde{\overline{u}}_{i,j+1}-\tilde{\overline{u}}_{i,j})^2,\\
\beta_{12}&=(\tilde{\overline{u}}_{i+1,j}-\tilde{\overline{u}}_{i,j})^2
           +(\tilde{\overline{u}}_{i,j+1}-\tilde{\overline{u}}_{i,j})^2,\\
\beta_{13}&=(\tilde{\overline{u}}_{i+1,j}-\tilde{\overline{u}}_{i,j})^2
           +(\tilde{\overline{u}}_{i,j}-\tilde{\overline{u}}_{i,j-1})^2,\\
\beta_{14}&=(\tilde{\overline{u}}_{i,j}-\tilde{\overline{u}}_{i-1,j})^2
           +(\tilde{\overline{u}}_{i,j}-\tilde{\overline{u}}_{i,j-1})^2,
\end{aligned}
\right.
\end{align}
and the absolute difference $\tau_{1}$ among these smoothness indicators
\begin{equation}
\label{2D-tau1}
\tau_1=\left(\frac{\sum_{k\neq l}|\beta_{1k}-\beta_{1l}|}{6}\right)^2,
\end{equation}
where the selection of the power is to be consistent with the definition (\ref{tau1}) of $\tau_{1}$ in the one dimension.
\newline
(2) We give these four polynomials $q_{1k}(x,y)$ the same linear weights $\gamma_{1k}=\frac{1}{4}$ for $k=1,2,3,4$ and calculate the corresponding nonlinear weights as
\begin{equation}
\label{2D-w1-1234}
\omega_{1k}=\frac{\overline{\omega}_{1k}}{\sum_{l=1}^{4}\overline{\omega}_{1l}},
\end{equation}
\begin{equation}
\label{2D-w1_-1234}
\overline{\omega}_{1k}=\gamma_{1k}\left(1+\frac{\tau_1}{\beta_{1k}+\varepsilon}\right),
\quad k=1,2,3,4,
\end{equation}
where $\varepsilon$ is still taken to be $10^{-10}$ as in the one dimensional case.
\newline
(3) We obtain a new polynomial
\begin{equation}
\label{2D-q1_new}
q_1^{new}(x,y)=\sum_{l=1}^{4}\omega_{1l}q_{1l}(x,y),
\end{equation}
and set $\beta_1$ to be
\begin{equation}
\label{2D-beta1}
\beta_{1}=\sum_{|\alpha|=1}|I_{i,j}|^{|\alpha|}
\left(\frac{\partial^{|\alpha|}}{\partial x^{{\alpha}_x}\partial y^{{\alpha}_y}}
q_1^{new}(x,y)\right)^2,
\end{equation}
where $\alpha=(\alpha_x,\alpha_y)$, $|\alpha|=\alpha_x+\alpha_y$.

Then, we define the corresponding nonlinear weights as (\ref{non-linear-weight}).

{\bf{\em Step 1.3.}} Obtain an approximation polynomial $u_{i,j}(x,y)$ of $u(x,y)$.

The new reconstruction polynomial $u_{i,j}(x,y)$ of $u(x,y)$ is defined as
\begin{equation}
\label{2D-re-p}
u_{i,j}(x,y)=\sum_{l=1}^{4}\omega_{l,4}p_l(x,y),
\end{equation}
and the Gauss-Lobatto point values that we need are taken to be
\begin{align}
\left.
\begin{aligned}
\label{2D-re-u}
u_{i \mp 1/2,j+\eta_{l}}^{\pm}&=u_{i,j}(x_{i \mp 1/2},y_{j+\eta_{l}}), \,\,\, l=1,2,3,4; \\
u_{i+\xi_{k},j \mp  1/2}^{\pm}&=u_{i,j}(x_{i+\xi_{k}},y_{j \mp  1/2}), \,\,\, k=1,2,3,4; \\
u_{i+\xi_{k},j+\eta_{l}}      &=u_{i,j}(x_{i+\xi_{k}},y_{j+\eta_{l}}), \,\,\, k=2,3; \,\, l=2,3,
\end{aligned}
\right.
\end{align}
where $\xi_{1}=\eta_{1}=-1/2$, $\xi_{2}=\eta_{2}=-\sqrt{5}/10$, $\xi_{3}=\eta_{3}=\sqrt{5}/10$ and $\xi_{4}=\eta_{4}=1/2$.

{\bf{\em Step 2.}} Update the Gauss-Lobatto point values of $u$ in the troubled-cells.

Here, we still choose the KXRCF troubled-cell indicator to identify the troubled-cells as in the one dimensional case, that is the target cell $I_{i,j}$ is identified to be a troubled-cell, if
\begin{equation}
\label{KXRCF_}
\aleph_{i,j}=\frac{\big|\int_{\partial I_{i,j}^{-}}\big(u_{i,j}(x,y)-u_{n_{i,j}}(x,y)\big)ds\big|}
{h_{i,j}^{\frac{l+1}{3}}\big|\partial I_{i,j}^{-}\big| ||u_{i,j}(x,y)||}>1,
\end{equation}
where $\partial I_{i,j}^{-}$ is the inflow boundary ($\overrightarrow{v}\cdot \overrightarrow{n}<0$, $\overrightarrow{v}$ is the velocity of the flow and $\overrightarrow{n}$ is the outer normal vector to $\partial I_{i,j}$), $I_{n_{i,j}}$ is the neighbor of $I_{i,j}$ on the side of $\partial I_{i,j}^{-}$, $h_{i,j}$ is the length of the cell $I_{i,j}$, the parameter $l$
(i.e. the degree of $u_{i,j}(x,y)$) is also taken to be 5, $u_{i,j}(x,y)$ is the approximation polynomial of $u(x,y)$ obtained in {\bf{\em Step 1}} and the norm is still
taken to be the $L^{\infty}$ norm.
If the target cell $I_{i,j}$ is identified to be a troubled-cell, we would like to modify the first order moment $\tilde{\overline{v}}_{i,j}$ in the $x$ direction by using the information of
$\{\tilde{\overline{u}}_{i-1,j},
   \tilde{\overline{u}}_{i  ,j},
   \tilde{\overline{u}}_{i+1,j},
   \tilde{\overline{v}}_{i-1,j},
   \tilde{\overline{v}}_{i+1,j}\}$
and modify the first order moment $\tilde{\overline{w}}_{i,j}$ in the $y$ direction  by using the information of
$\{\tilde{\overline{u}}_{i,j-1},
   \tilde{\overline{u}}_{i,j  },
   \tilde{\overline{u}}_{i,j+1},
   \tilde{\overline{w}}_{i,j-1},
   \tilde{\overline{w}}_{i,j+1}\}$ in a dimension-by-dimension manner.
After modifying the first order moments in the troubled-cells, we repeat the reconstruction process {\bf{\em Step 1}}
for these troubled-cells to update the corresponding Gauss-Lobatto point values of $u$.

{\bf{\em Step 3.}} Discretize the semi-discrete scheme in time.

After all these Gauss-Lobatto point values are obtained, we put them into the formula of the numerical flux.
Then, we discretize (\ref{2D-equ-int-approximation}) by the third-order TVD Runge-Kutta
method (\ref{Runge-Kutta}) in time to complete the entire discretization process.

{\bf{\em Remark 3.}}
In {\bf{\em Step 2}} above, we still choose the KXRCF troubled-cell indicator to identify troubled-cells.
As shown in Section 3, the KXRCF troubled-cell indicator works pretty well for our scheme in the two-dimensional
case as well. What needs a special attention is that,
for the two dimensional scalar equation,
     the solution $u$ is defined as  our indicator variable, and then the corresponding
     $\overrightarrow{v}=f^{'}(u)$ in the $x$-direction and $\overrightarrow{v}=g^{'}(u)$ in the $y$-direction;
for the two dimensional Euler system,
only the density $\rho$ is set to be our indicator variable, and then the corresponding
     $\overrightarrow{v}=\mu$ is the velocity in the $x$-direction of the fluid and
$\overrightarrow{v}=\nu$ is the velocity in the $y$-direction of the fluid.
In short, for the two dimensional case, the line integral average in the formula (\ref{KXRCF_}) is approximated by a four-point Gauss-Lobatto integration, that is
\begin{align}
\label{line-integral-average_}
\frac{1}{\big|\partial I_{i,j}^{-}\big|}\left|\int_{\partial I_{i,j}^{-}}\big(u_{i,j}(x,y)-u_{n_{i,j}}(x,y)\big)ds\right| \notag
&=\Bigg|
  \sum_{l=1}^4\triangle y\omega_l
                      (u_{i-\frac{1}{2},j+\eta_{l}}^{+}
                      -u_{i-\frac{1}{2},j+\eta_{l}}^{-})
                      *sf(\overrightarrow{v}_{i-\frac{1}{2},j})\\ \notag
&+\sum_{l=1}^4\triangle y\omega_l
                      (u_{i+\frac{1}{2},j+\eta_{l}}^{-}
                      -u_{i+\frac{1}{2},j+\eta_{l}}^{+})
                      *sf(-\overrightarrow{v}_{i+\frac{1}{2},j})\\ \notag
&+\sum_{k=1}^4\triangle x\omega_k
                      (u_{i+\xi_{k},j-\frac{1}{2} }^{+}
                      -u_{i+\xi_{k},j-\frac{1}{2} }^{-})
                      *sf(\overrightarrow{v}_{i,j-\frac{1}{2}})\\
&+\sum_{k=1}^4\triangle x\omega_k
                      (u_{i+\xi_{k},j+\frac{1}{2} }^{-}
                      -u_{i+\xi_{k},j+\frac{1}{2} }^{+})
                      *sf(-\overrightarrow{v}_{i,j+\frac{1}{2}})
\Bigg|\frac{1}{\big|\partial I_{i,j}^{-}\big|},
\end{align}

and the norm $||u_{i,j}(x,y)||$ is
taken to be the maximum norm of all the Gauss-Lobatto point values in the cell $I_{i,j}
$(i.e. $||u_{i,j}(x,y)||\approx \max$$\{ |u_{i\mp1/2,j+\eta_{l}}^{\pm}|: \, l=1,2,3,4$;
$|u_{i+\xi_{k},j\mp1/2}^{\pm}|: \, k=1,2,3,4$;
$|u_{i+\xi_{k},j+\eta_{l}}|: \, k=2,3;l=2,3\}$).
Note that all the values used in the troubled-cell indicator are also already obtained in the reconstruction process {\bf{\em
Step 1}}, thus there is no need to reconstruct an extra polynomial as in \cite{zcq2020,zq2020}.
}

\section{Numerical tests}
\label{sec3}
\setcounter{equation}{0}
\setcounter{figure}{0}
\setcounter{table}{0}

In this section, a number of typical numerical examples are
given to demonstrate the stability and resolution of our moment-based multi-resolution HWENO scheme which is termed as ``HWENO6-M5-I/NI" where ``HWENO6" means the reconstruction process is of the sixth order of accuracy, ``M5" means the first order moments of the troubled-cells are modified by a quartic polynomial, and ``I" means only the first order moments of the troubled-cells are modified while ``NI" means the first order moments of all the cells are modified without judgment.
Before we start to show the results of the examples we have calculated, let us first explain some of the parameters
in particular: the first one is that we set the CFL number as 0.6 for both the one and two dimensional cases, but note that
for the accuracy tests a suitably reduced time step is used in order to ensure the dominance of the spatial error; the
second one is that we take the linear weights as
$\overline{\gamma}_{1,4}=1  $,     $\overline{\gamma}_{2,4}=10  $, $\overline{\gamma}_{3,4}=100$  and $\overline{\gamma}_{4,4}=1000$
both in the one and two dimensions in this paper.

\bigskip

\begin{figure}[htbp]
  \centering{\includegraphics[width=2.5 in,height=2.5 in]{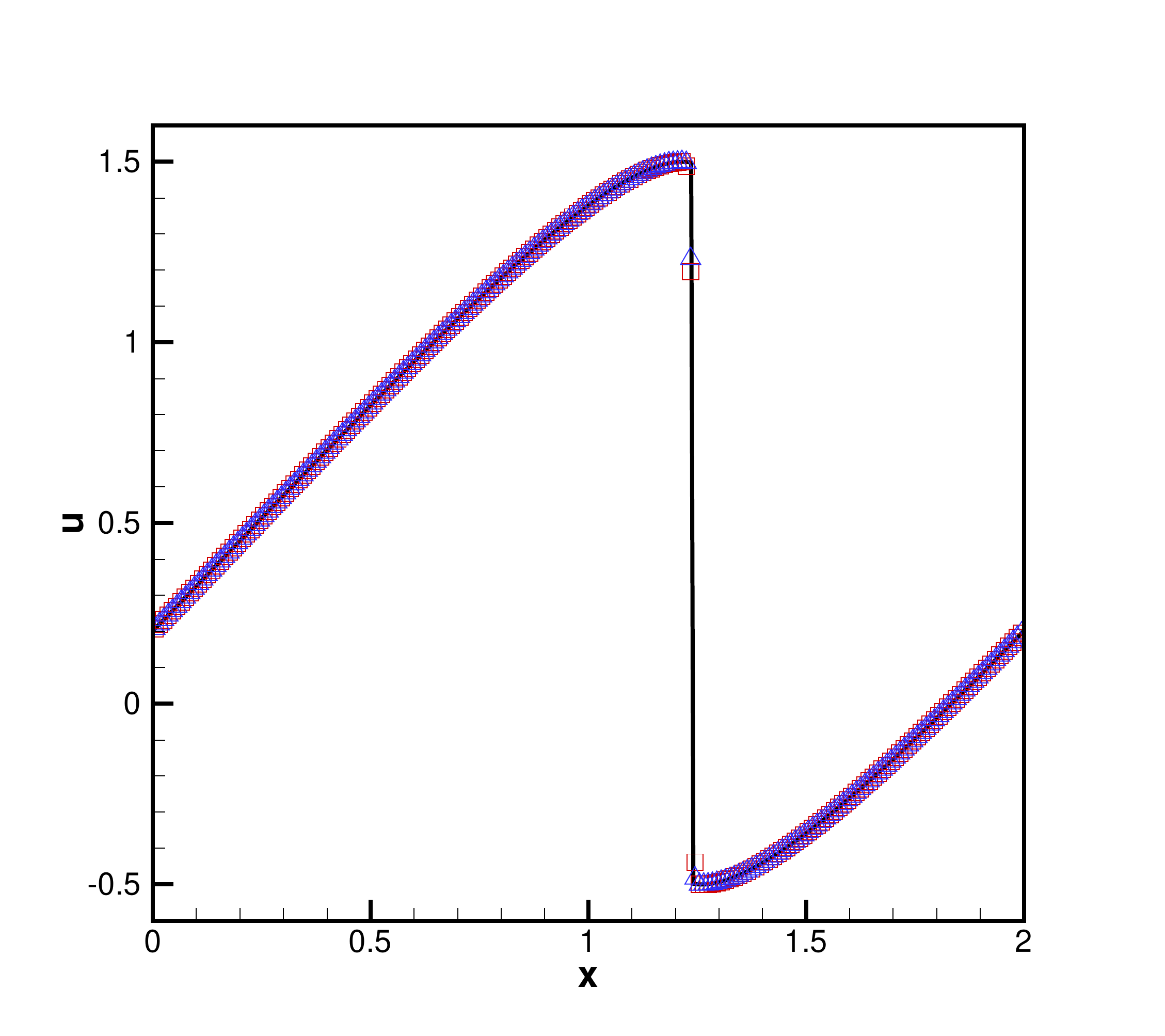}
             \includegraphics[width=2.5 in,height=2.5 in]{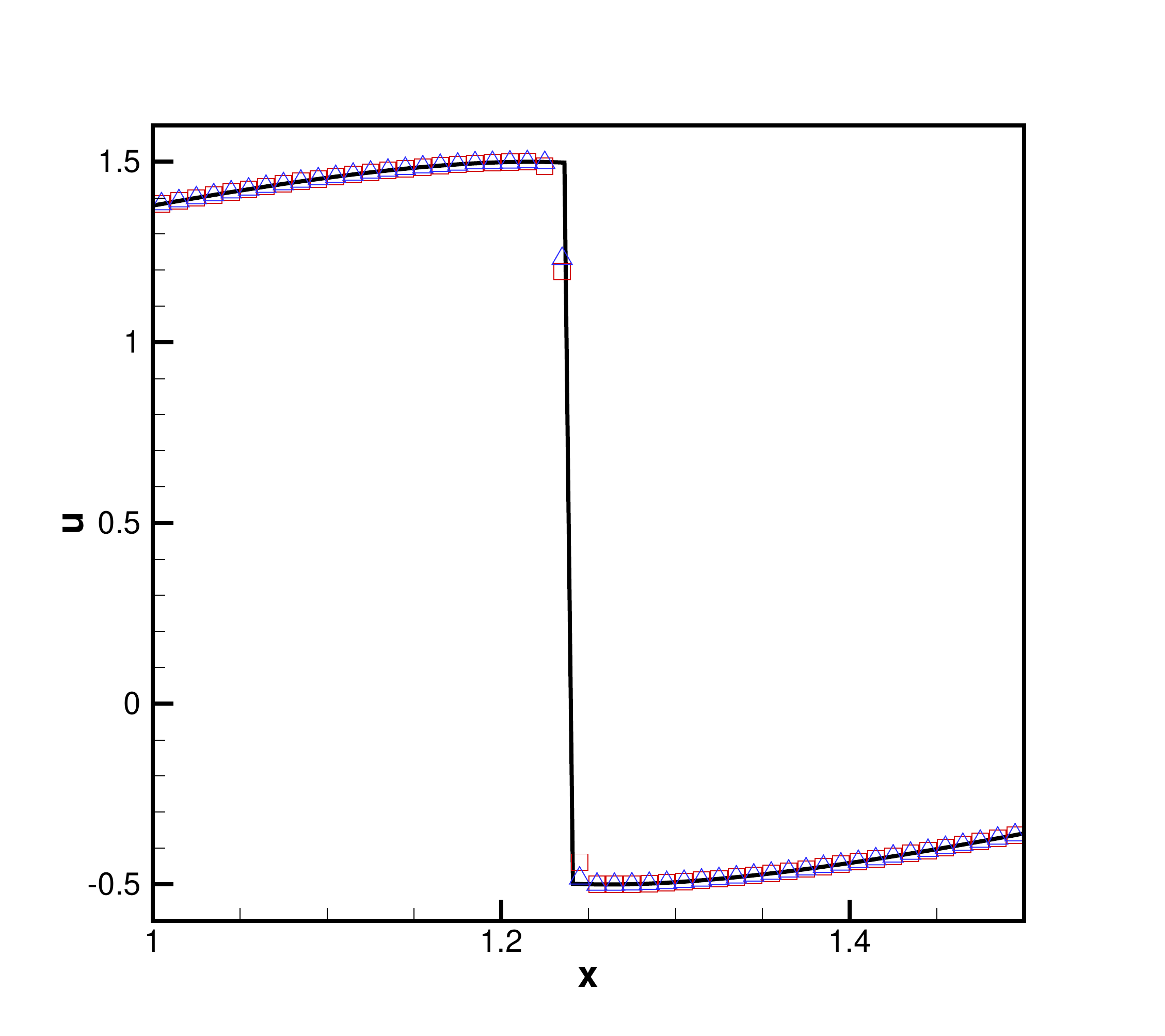}}
  \caption{{
  1D-Burgers' equation.
  $T=1.5/\pi$.
  Left:  density;
  right: density zoomed in.
  Solid line: the exact solution;
  squares:    the result of the HWENO6      scheme;
  triangles:  the result of the HWENO6-M5-I scheme.
  Number of cells: 200.}
\label{ex1}}
\end{figure}

\noindent {\bf Example 3.1.}
One-dimensional scalar Burgers' equation:
\begin{equation}
\label{1D-Burgers}
\mu_t+\left(\frac{\mu^2}{2}\right)_x=0, \ \ 0<x<2,
\end{equation}
with the initial condition $\mu(x,0)=0.5+\sin(\pi x)$
and periodic boundary condition.
As we know, when the time reaches $T=0.5/\pi$, its solution is still smooth, and now the corresponding errors and orders obtained by the
HWENO6-M5-I and HWENO6-M5-NI schemes are listed in the Table \ref{table1};
but when the time reaches $T=1.5/\pi$, there will appear discontinuities in its solution, and now the corresponding numerical solutions
obtained by the HWENO6 and HWENO6-M5-I schemes are drawn on the Fig \ref{ex1} in comparison with the reference exact solution.

\begin{table}[ht]
\begin{center}
{\caption{
$\mu_t+\left(\frac{\mu^2}{2}\right)_x=0$.
The initial condition $\mu(x,0)=0.5+\sin(\pi x)$.
Periodic boundary condition.
$T=0.5/\pi$.
HWENO6-M5-I and HWENO6-M5-NI schemes.
$L^1$ and $L^{\infty}$.}
\begin{tabular} {|c|c|c|c|c|c|c|c|c|} \hline
 &\multicolumn{4}{|c|}{HWENO6-M5-I  scheme}
 &\multicolumn{4}{|c|}{HWENO6-M5-NI scheme}\\
      \hline
  grid points &$ L^1$ error &  order & $L^\infty$ error &  order &
               $ L^1$ error &  order & $L^\infty$ error &  order
   \\ \hline
    10&5.24E-03   &       &2.26E-02   &       &6.48E-03   &       &3.01E-02   &           \\ \hline
    20&5.67E-04   &3.21   &4.61E-03   &2.29   &3.11E-04   &4.38   &9.84E-04   &4.93       \\ \hline
    40&1.23E-06   &8.85   &1.37E-05   &8.40   &6.97E-06   &5.48   &2.55E-05   &5.27       \\ \hline
    80&2.90E-09   &8.72   &3.27E-08   &8.71   &5.96E-08   &6.87   &9.36E-07   &4.77       \\ \hline
   160&2.46E-11   &6.88   &3.52E-10   &6.54   &1.88E-09   &4.98   &2.95E-08   &4.99       \\ \hline
   320&3.95E-13   &5.96   &3.69E-12   &6.57   &6.02E-11   &4.97   &9.53E-10   &4.95       \\ \hline
\end{tabular}
\label{table1}}
\end{center}
\end{table}

\smallskip

\bigskip

\begin{figure}[htbp]
  \centering{\includegraphics[width=2.5 in,height=2.5 in]{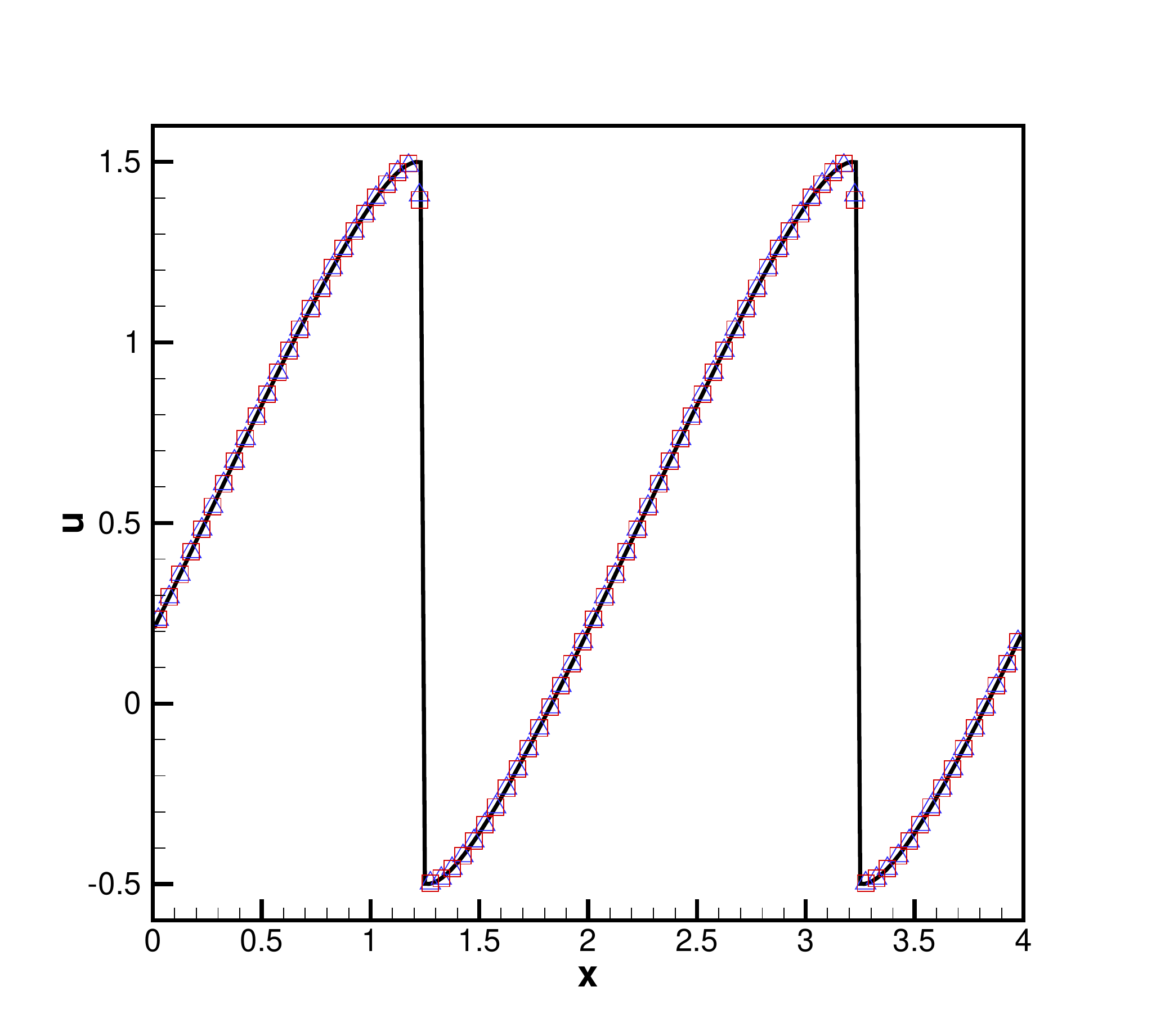}
             \includegraphics[width=2.5 in,height=2.5 in]{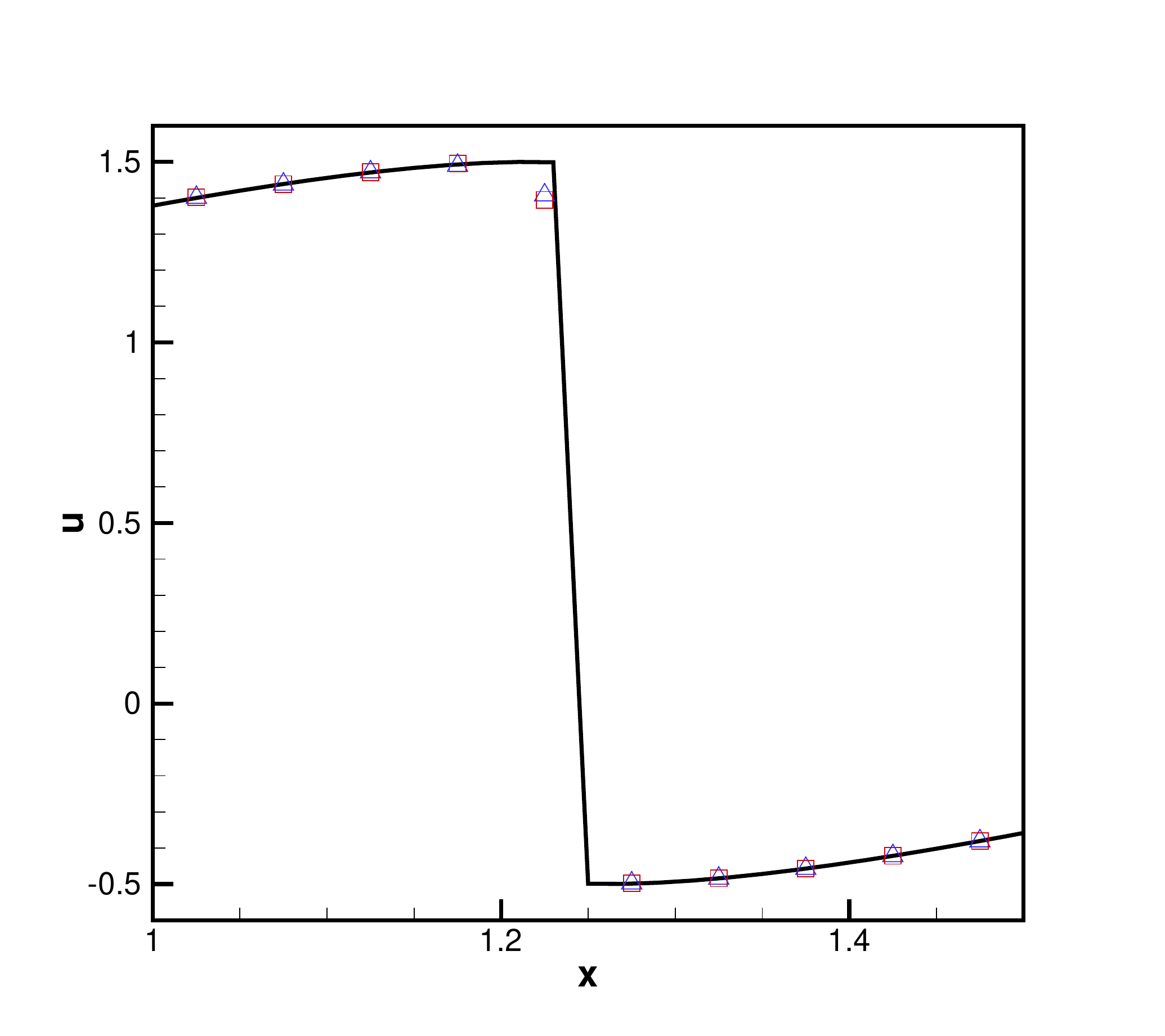}}
  \caption{{
  2D-Burgers' equation.
  $T=1.5/\pi$.
  Left:  density;
  right: density zoomed in at $x=y$.
  Solid line: the exact solution;
  squares:    the result of the HWENO6      scheme;
  triangles:  the result of the HWENO6-M5-I scheme.
  Number of cells: 200$\times$200.}
\label{ex2}}
\end{figure}

\noindent {\bf Example 3.2.}
Two-dimensional scalar Burgers' equation:
\begin{equation}
\label{2D-Burgers} \mu_t+\left(\frac{\mu^2}{2}\right)_x+\left(\frac{\mu^2}{2}\right)_y=0, \ \ 0<x,y<4,
\end{equation}
with the initial condition $\mu(x,y,0)=0.5+\sin\big(\pi (x+y)/2\big)$
and periodic boundary condition.
As is known to us, when the time reaches $T=0.5/\pi$, its solution is still smooth, and now the corresponding errors and orders obtained by the
HWENO6-M5-I and HWENO6-M5-NI schemes are listed in the Table \ref{table2}; but when the time reaches $T=1.5/\pi$, there will appear discontinuities in its solution, and now the corresponding numerical solutions at $x=y$ obtained by the HWENO6 and HWENO6-M5-I schemes are drawn on the Fig \ref{ex2} in comparison with the reference exact solution.

\begin{table}[ht]
\begin{center}
{\caption{
$\mu_t+\left(\frac{\mu^2}{2}\right)_x+\left(\frac{\mu^2}{2}\right)_y=0$.
The initial condition $\mu(x,y,0)=0.5+\sin\big(\pi (x+y)/2\big)$.
Periodic boundary condition.
$T=0.5/\pi$.
HWENO6-M5-I and HWENO6-M5-NI schemes.
$L^1$ and $L^{\infty}$.}
\begin{tabular} {|c|c|c|c|c|c|c|c|c|} \hline
 &\multicolumn{4}{|c|}{HWENO6-M5-I  scheme}
 &\multicolumn{4}{|c|}{HWENO6-M5-NI scheme}\\
      \hline
  grid points &$ L^1$ error &  order & $L^\infty$ error &  order &
               $ L^1$ error &  order & $L^\infty$ error &  order
   \\ \hline
    10$\times$ 10 &1.03E-02 &     &2.96E-02 &     &1.08E-02 &     &4.31E-02 &           \\ \hline
    20$\times$ 20 &7.78E-05 &7.05 &2.65E-04 &6.80 &1.42E-04 &6.25 &8.03E-04 &5.74       \\ \hline
    40$\times$ 40 &2.41E-06 &5.01 &1.56E-05 &4.09 &6.62E-06 &4.42 &3.84E-05 &4.39       \\ \hline
    80$\times$ 80 &3.86E-08 &5.97 &5.42E-07 &4.84 &1.03E-07 &6.00 &8.33E-07 &5.53       \\ \hline
   160$\times$160 &5.62E-10 &6.10 &9.83E-09 &5.78 &2.87E-09 &5.17 &2.30E-08 &5.18       \\ \hline
   320$\times$320 &6.83E-12 &6.36 &1.24E-10 &6.31 &9.77E-11 &4.88 &7.42E-10 &4.95       \\ \hline
\end{tabular}
\label{table2}}
\end{center}
\end{table}

\smallskip

\bigskip

\noindent {\bf Example 3.3.}
One-dimensional Euler equations:
\begin{equation}
\label{1D-Euler}
\frac{\partial}{\partial t}
\left(
\begin{array}{c}
 \rho \\ \rho\mu \\ E
\end{array}
\right)+
\frac{\partial}{\partial x}
\left(
\begin{array}{c}
 \rho\mu \\ \rho\mu^2+p \\ \mu(E+p)
\end{array}
\right)=0, \ \ 0<x<2\pi,
\end{equation}
where $\rho$ is the density, $\mu$ is the velocity, $E$ is the total energy and $p$ is the pressure.
The initial conditions are
\begin{equation}
\label{1D-initial-condition}
\begin{array}{l}
\rho(x,0)=\frac{1+0.2\sin(x)}{2\sqrt3}, \,\,
\mu(x,0)=\sqrt\gamma\rho(x,0), \,\,
p(x,0)=\rho(x,0)^{\gamma},
\end{array}
\end{equation}
and the boundary conditions are periodic.
The exact solution of above Euler equations is given in \cite{lqs2021}.
When $T=3$, its solution is still smooth, and now the corresponding errors and orders obtained by the
HWENO6-M5-I and HWENO6-M5-NI schemes are listed in the Table \ref{table3}.

\begin{table}[ht]
\begin{center}
{\caption{
1D-Euler equations:
The initial condition $\rho(x,0)=\frac{1+0.2\sin(x)}{2\sqrt3}$, $\mu(x,0)=\sqrt\gamma\rho(x,0)$ and $p(x,0)=\rho(x,0)^{\gamma}$.
Periodic boundary condition.
$T=3$.
HWENO6-M5-I and HWENO6-M5-NI schemes.
$L^1$ and $L^{\infty}$.}
\begin{tabular} {|c|c|c|c|c|c|c|c|c|} \hline
 &\multicolumn{4}{|c|}{HWENO6-M5-I  scheme}
 &\multicolumn{4}{|c|}{HWENO6-M5-NI scheme}\\
      \hline
  grid points &$ L^1$ error &  order & $L^\infty$ error &  order &
               $ L^1$ error &  order & $L^\infty$ error &  order
   \\ \hline
   10&8.54E-04   &       &2.95E-03   &       &8.76E-04   &       &3.67E-03   &           \\ \hline
   20&1.87E-05   &5.52   &1.41E-04   &4.39   &1.09E-04   &3.01   &5.35E-04   &2.78       \\ \hline
   40&4.99E-07   &5.23   &9.01E-06   &3.97   &8.44E-06   &3.69   &1.05E-04   &2.35       \\ \hline
   80&7.57E-09   &6.04   &1.95E-07   &5.53   &2.49E-07   &5.08   &5.37E-06   &4.28       \\ \hline
  160&1.10E-10   &6.10   &2.60E-09   &6.23   &7.12E-09   &5.13   &1.66E-07   &5.02       \\ \hline
  320&1.69E-12   &6.03   &3.65E-11   &6.15   &2.15E-10   &5.05   &5.03E-09   &5.05       \\ \hline
\end{tabular}
\label{table3}}
\end{center}
\end{table}
\smallskip

\bigskip

\noindent {\bf Example 3.4.}
Two-dimensional Euler equations:
\begin{equation}
\label{2D-Euler}
\frac{\partial}{\partial t}
\left(
\begin{array}{c}
 \rho \\ \rho\mu \\ \rho\nu \\ E
\end{array}
\right)+
\frac{\partial}{\partial x}
\left(
\begin{array}{c}
 \rho\mu \\ \rho\mu^2+p \\ \rho\mu\nu \\ \mu(E+p)
\end{array}
\right)+
\frac{\partial}{\partial y}
\left(
\begin{array}{c}
 \rho\nu \\ \rho\mu\nu \\ \rho\nu^2+p \\ \nu(E+p)
\end{array}
\right)=0, \ \ 0<x,y<4\pi,
\end{equation}
where $\rho$ is the density, $\mu$ is the velocity in the $x$-direction, $\nu$ is the velocity in the $y$-direction, $E$ is the total energy and $p$ is the pressure.
The initial conditions are
\begin{equation}
\label{2D-initial-condition}
\begin{array}{l}
\rho(x,y,0)=\frac{1+0.2\sin(\frac{x+y}{2})}{\sqrt6}, \,\,
\mu(x,y,0)=\nu(x,y,0)=\sqrt{\frac{\gamma}{2}}\rho(x,y,0), \,\,
p(x,y,0)=\rho(x,y,0)^{\gamma},
\end{array}
\end{equation}
and the boundary conditions are periodic in both directions.
The exact solution of above Euler equations is given in \cite{lqs2021}.
When $T=3$, its solution is still smooth, and now the corresponding errors and orders obtained by the
HWENO6-M5-I and HWENO6-M5-NI schemes are listed in the Table \ref{table4}.

\begin{table}[ht]
\begin{center}
{\caption{
2D-Euler equations:
The initial condition
$\rho(x,y,0)=\frac{1+0.2\sin(\frac{x+y}{2})}{\sqrt6}$, $\mu(x,y,0)=\nu(x,y,0)=\sqrt{\frac{\gamma}{2}}\rho(x,y,0)$ and $p(x,y,0)=\rho(x,y,0)^{\gamma}$.
Periodic boundary condition.
$T=3$.
HWENO6-M5-I and HWENO6-M5-NI schemes.
$L^1$ and $L^{\infty}$.}
\begin{tabular} {|c|c|c|c|c|c|c|c|c|} \hline
 &\multicolumn{4}{|c|}{HWENO6-M5-I  scheme}
 &\multicolumn{4}{|c|}{HWENO6-M5-NI scheme}\\
      \hline
  grid points &$ L^1$ error &  order & $L^\infty$ error &  order &
               $ L^1$ error &  order & $L^\infty$ error &  order
   \\ \hline
    10$\times$ 10 &6.56E-03 &     &1.74E-02 &     &6.70E-03 &     &2.20E-02 &           \\ \hline
    20$\times$ 20 &2.54E-04 &4.69 &1.78E-03 &3.29 &4.20E-04 &4.00 &2.73E-03 &3.01       \\ \hline
    40$\times$ 40 &1.13E-05 &4.50 &1.13E-04 &3.97 &2.47E-05 &4.09 &2.23E-04 &3.61       \\ \hline
    80$\times$ 80 &2.70E-07 &5.38 &6.06E-06 &4.22 &7.57E-07 &5.03 &1.47E-05 &3.93       \\ \hline
   160$\times$160 &4.18E-09 &6.01 &1.14E-07 &5.73 &2.14E-08 &5.14 &4.76E-07 &4.95       \\ \hline
   320$\times$320 &4.67E-11 &6.49 &1.37E-09 &6.39 &6.24E-10 &5.10 &1.30E-08 &5.20       \\ \hline
\end{tabular}
\label{table4}}
\end{center}
\end{table}
\smallskip

\bigskip

\noindent {\bf Comment:}
According to the results listed in above four tables
Table \ref{table1}, Table \ref{table2}, Table \ref{table3} and Table \ref{table4},
we can see that the HWENO6-M5-I scheme can reach sixth order of accuracy, this is because the reconstruction process is of the sixth order and the modification procedure has not been enacted for these smooth cases.
If we modify the first order moments of all the cells without judgment, we can see that the HWENO6-M5-NI scheme can reach fifth order of accuracy as expected, this is because the first order moments of all the cells are modified by a quartic polynomial, which leads to the decrease in the order of accuracy.

From Fig \ref{ex1} and Fig \ref{ex2},
we can observe that both the HWENO6 and HWENO6-M5-I schemes work well in comparison with the exact solutions and there is not much difference between the results of these two schemes.
\smallskip

\bigskip

\begin{figure}[htbp]
  \centering{\includegraphics[width=2.5 in,height=2.5 in]{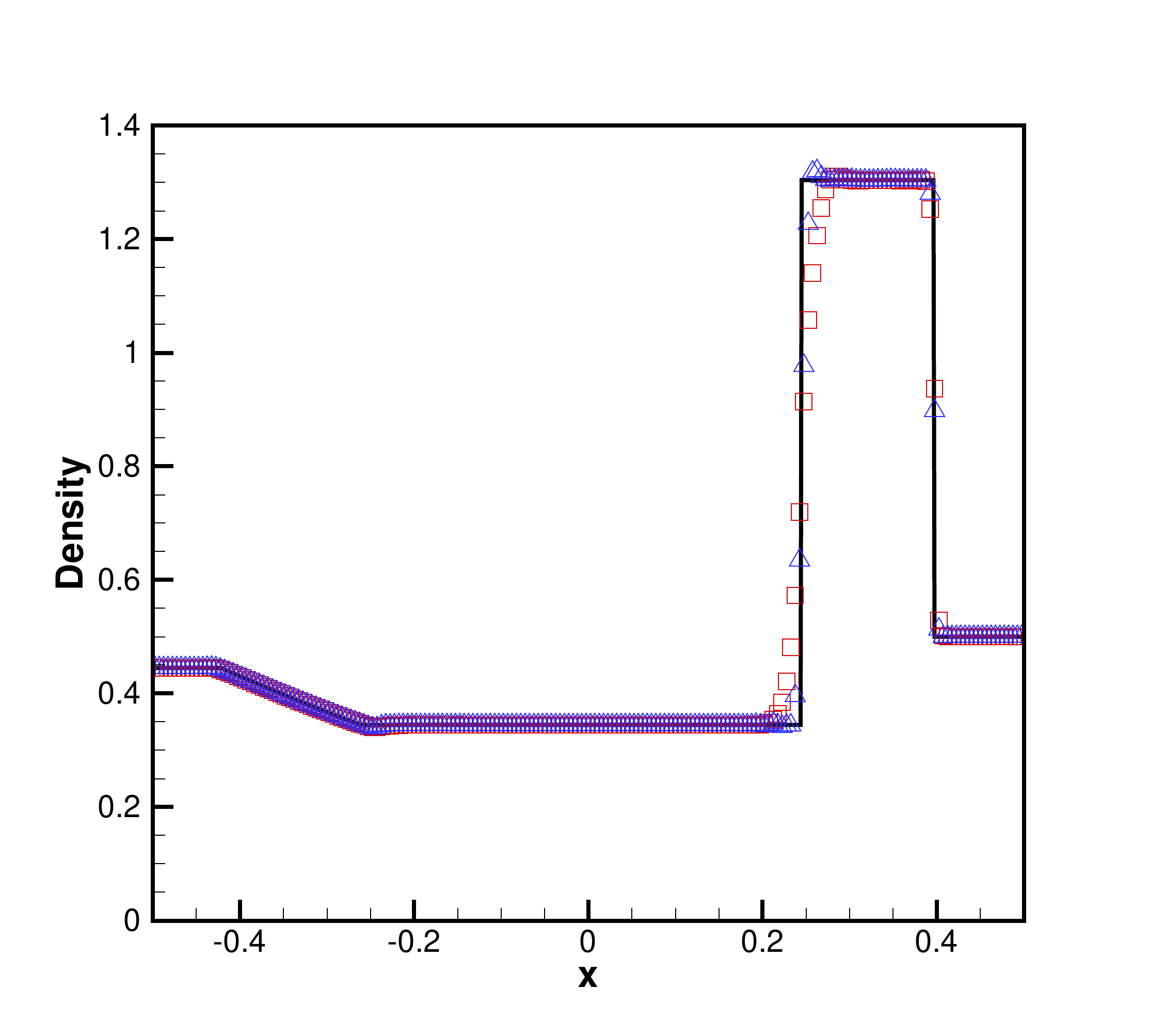}
             \includegraphics[width=2.5 in,height=2.5 in]{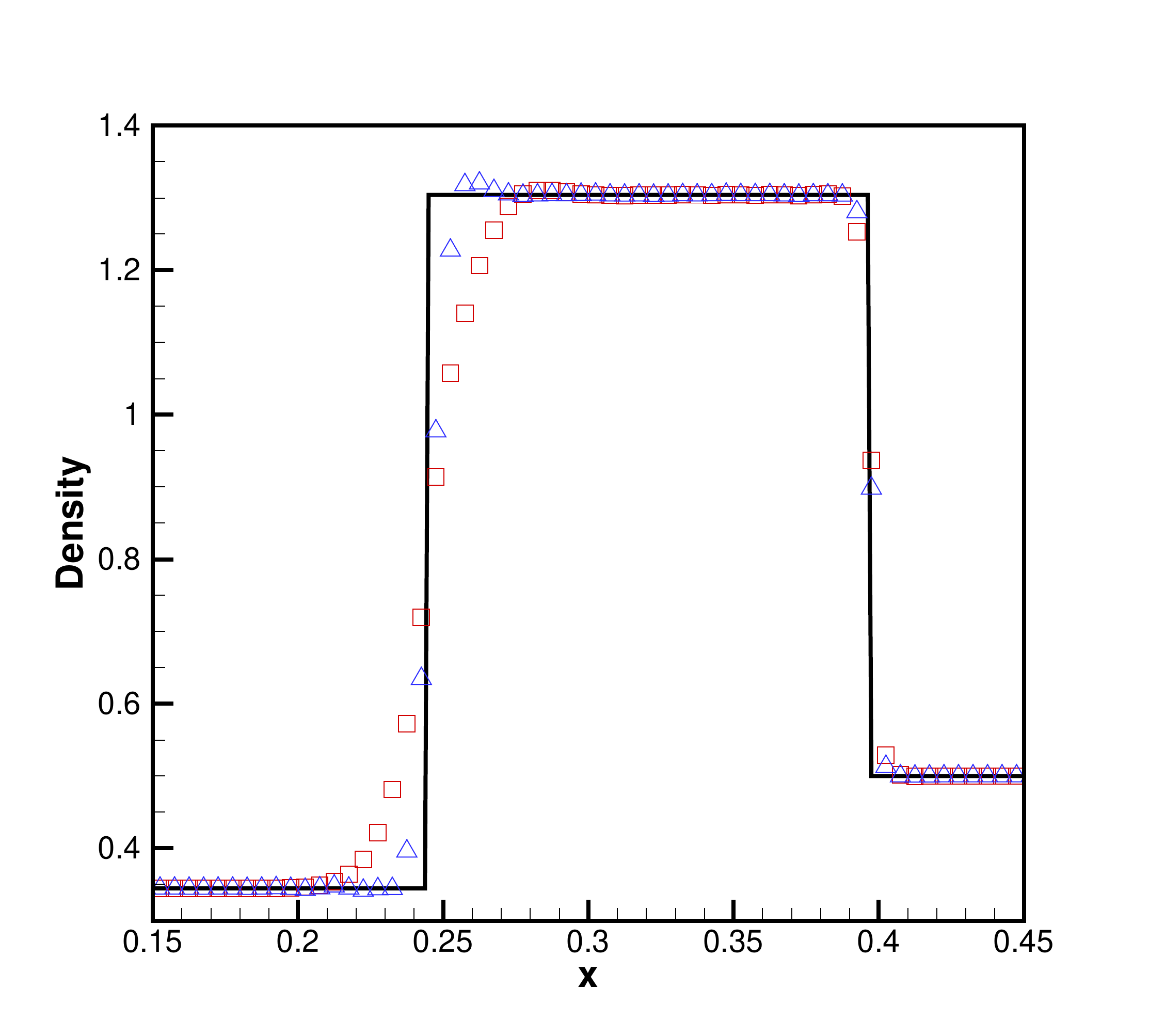}
             \includegraphics[width=2.5 in,height=2.5 in]{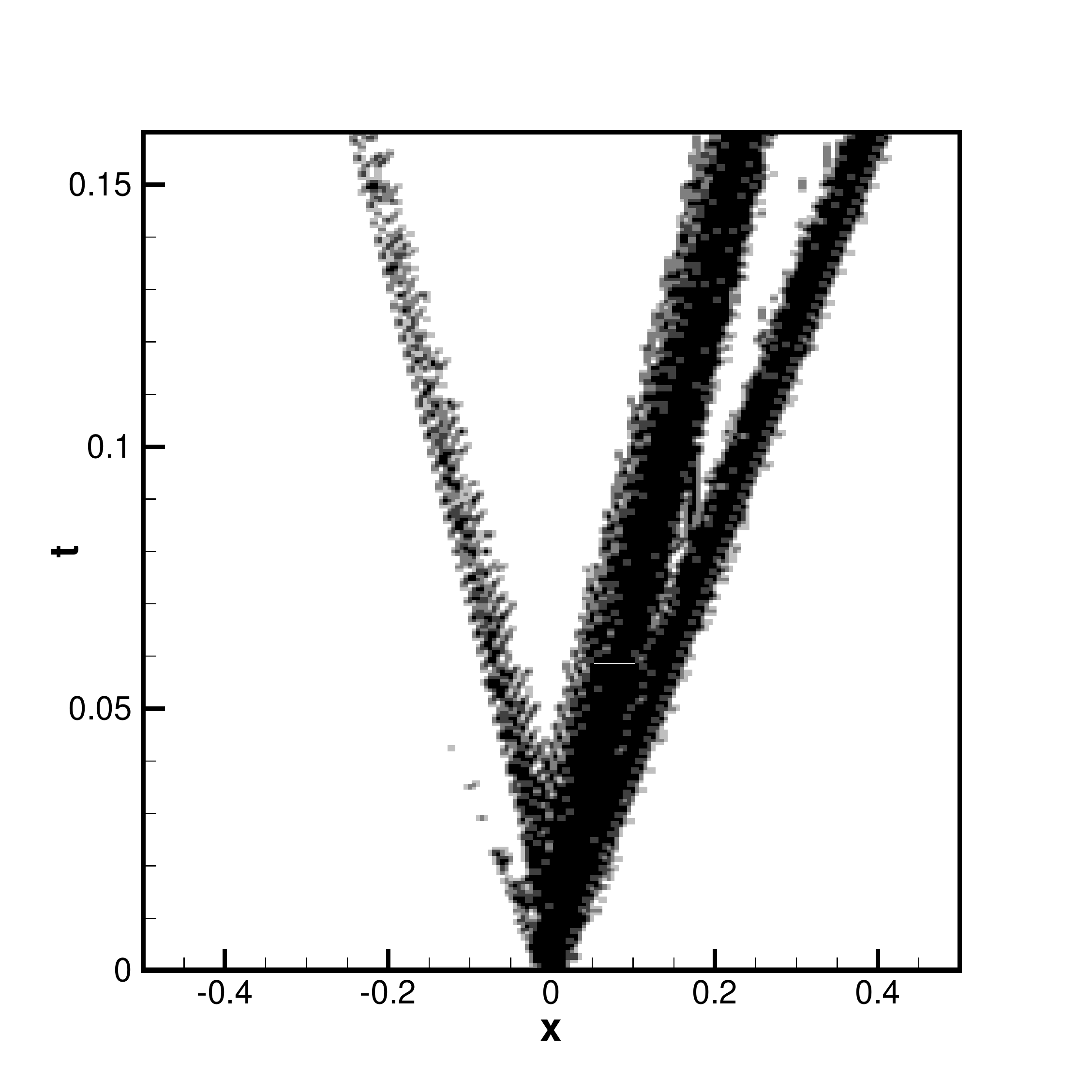}}
  \caption{{
  The Lax problem.
  $T=0.16$.
  Top:    density, density zoomed in;
  bottom: the locations of the troubled-cells.
  Solid line: the exact solution;
  squares:    the result of the HWENO6      scheme;
  triangles:  the result of the HWENO6-M5-I scheme.
  Number of cells: 200.}
\label{ex5}}
\end{figure}

\noindent {\bf Example 3.5.}
The Lax problem: one-dimensional Euler equations (\ref{1D-Euler}) with the Riemann initial condition:
\begin{equation}
\label{1D-Lax}
(\rho,\mu,p)^T=
\left\{
\begin{array}{c}
 (0.445,0.698,3.528)^T, \ \ -0.5<x<0  ,\\
 (0.5  ,0    ,0.571)^T, \ \  0  <x<0.5.
\end{array}
\right.
\end{equation}
The computed result of the density $\rho$ and its zoom-in picture obtained by the HWENO6 and HWENO6-M5-I schemes at the final time  $T=0.16$ in comparison with the reference exact solution, as well as the corresponding locations of the troubled-cells for the HWENO6-M5-I scheme over time are plotted in the Fig \ref{ex5}.
\smallskip

\bigskip

\begin{figure}[htbp]
  \centering{\includegraphics[width=2.5 in,height=2.5 in]{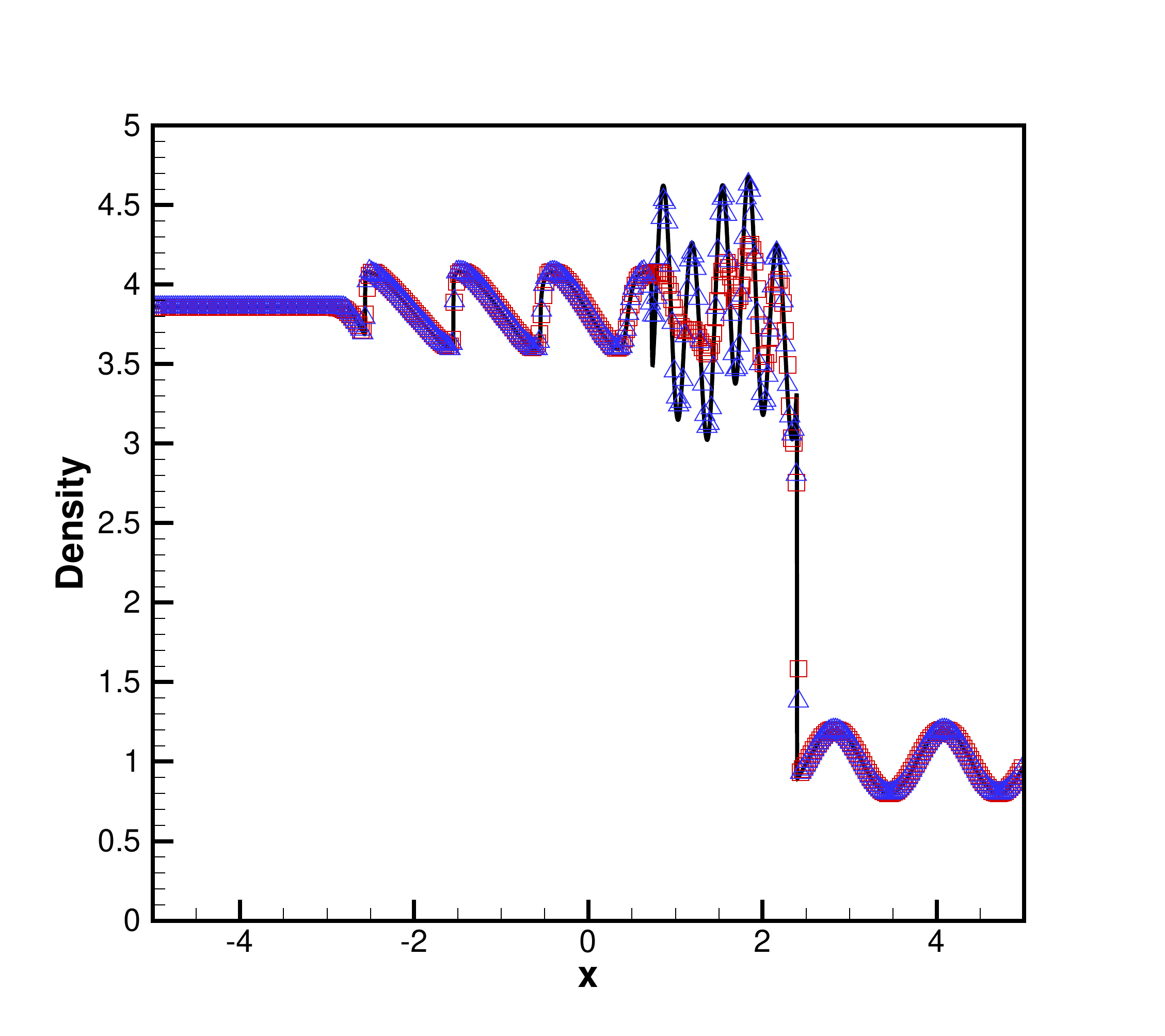}
             \includegraphics[width=2.5 in,height=2.5 in]{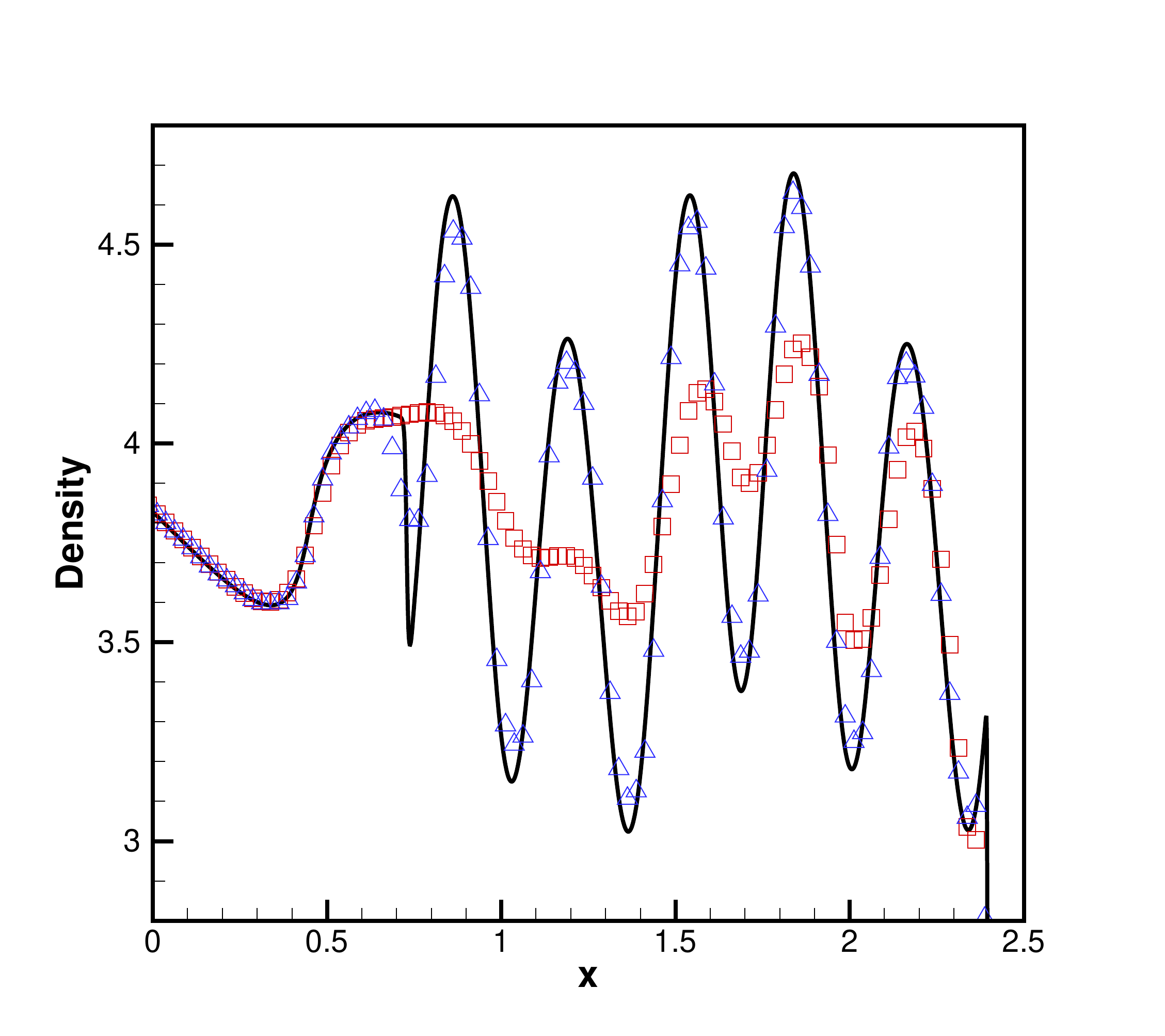}
             \includegraphics[width=2.5 in,height=2.5 in]{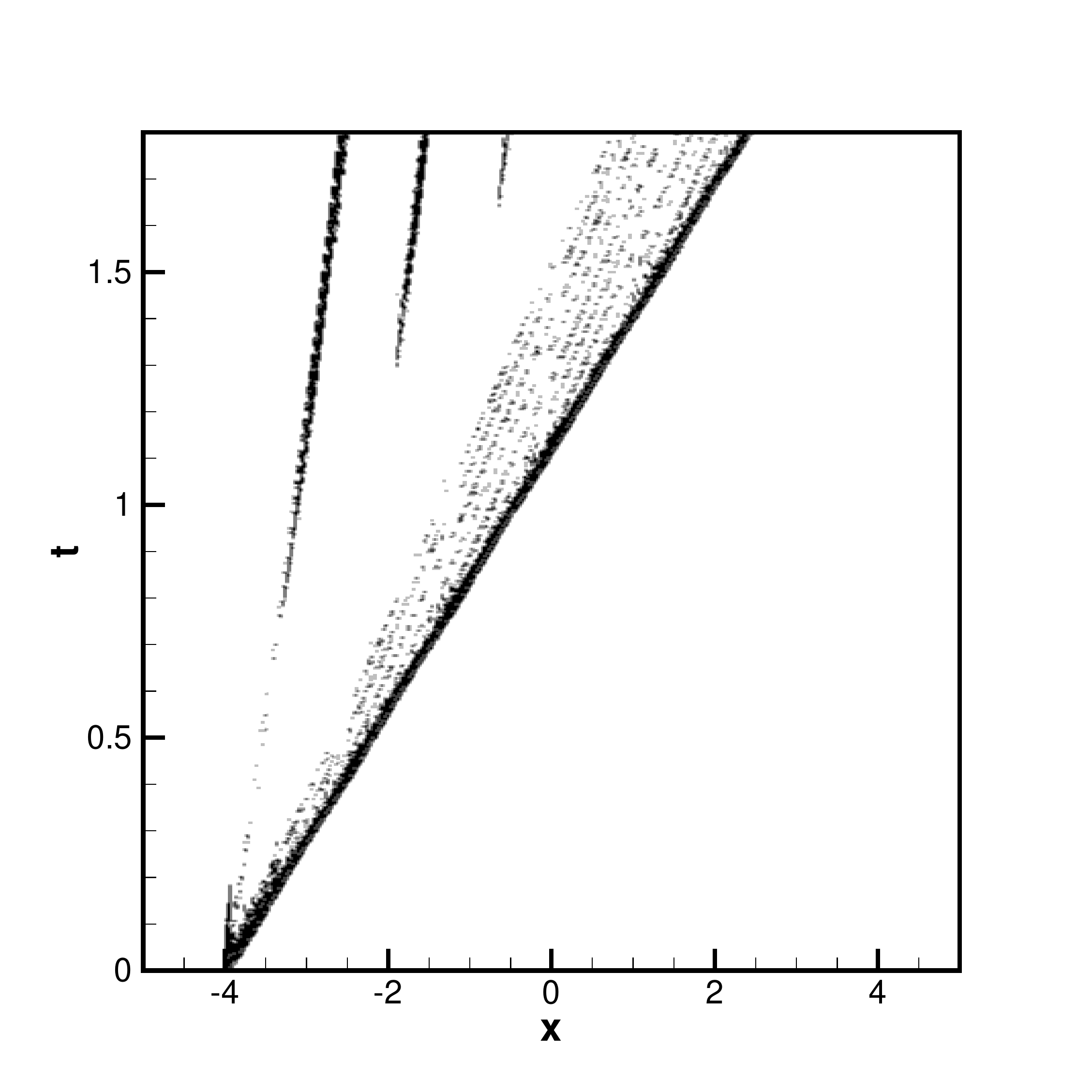}}
  \caption{{
  The shock density wave interaction problem.
  $T=1.8$.
  Top:    density, density zoomed in;
  bottom: the locations of the troubled-cells.
  Solid line: the exact solution;
  squares:    the result of the HWENO6      scheme;
  triangles:  the result of the HWENO6-M5-I scheme.
  Number of cells: 400.}
\label{ex6}}
\end{figure}

\noindent {\bf Example 3.6.}
The shock density wave interaction problem: one-dimensional Euler equations (\ref{1D-Euler}) with a moving Mach=3 shock interaction containing sine waves in the density:
\begin{equation}
\label{1D-Shock}
(\rho,\mu,p)^T=
\left\{
\begin{array}{c}
     (3.857143     ,2.629369,10.333333    )^T, \ \ -5<x<-4,\\
 \big(1+0.2\sin(5x),0       ,1        \big)^T, \ \ -4<x< 5.
\end{array}
\right.
\end{equation}
The computational result of the density $\rho$ and its zoom-in picture obtained by the HWENO6 and HWENO6-M5-I schemes at the final time $T=1.8$ in comparison with the reference ``exact" solution (which is a numerically
converged solution computed by the fifth-order finite difference WENO scheme \cite{js1996} with 8000 grid points), as well as the corresponding
locations of the troubled-cells for the HWENO6-M5-I scheme over time are plotted in the Fig \ref{ex6}.
\smallskip

\bigskip

\begin{figure}[htbp]
  \centering{\includegraphics[width=2.5 in,height=2.5 in]{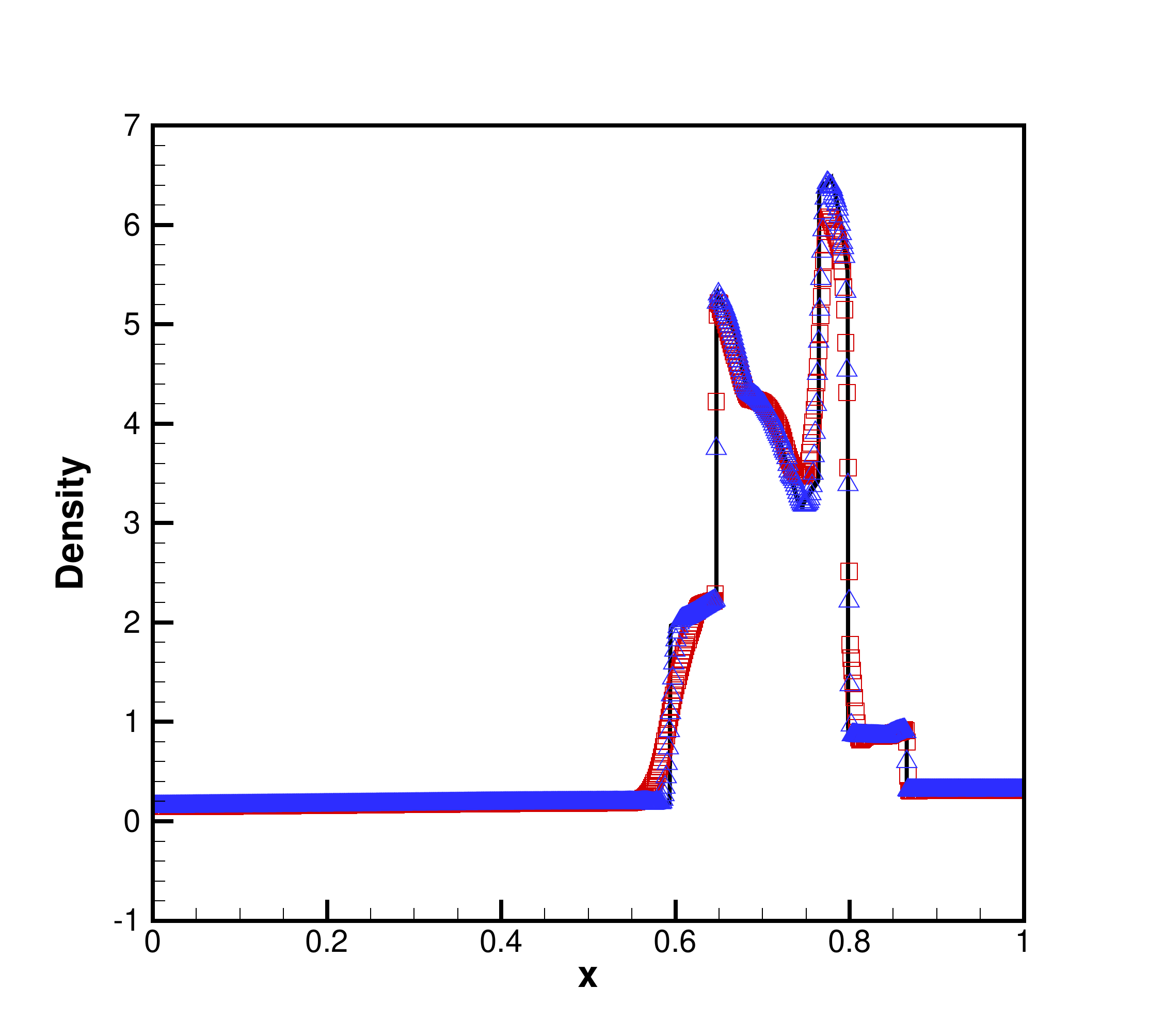}
             \includegraphics[width=2.5 in,height=2.5 in]{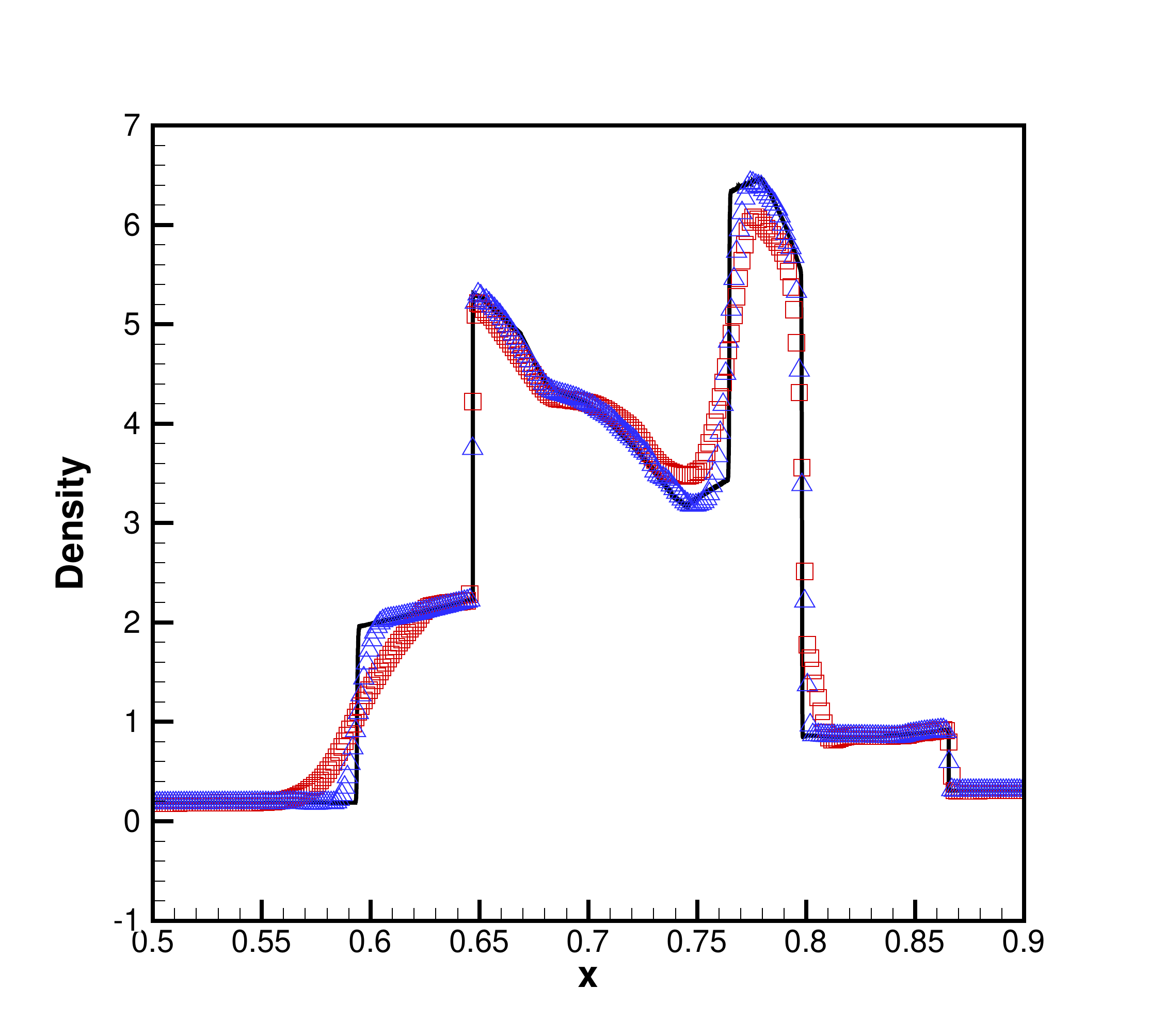}
             \includegraphics[width=2.5 in,height=2.5 in]{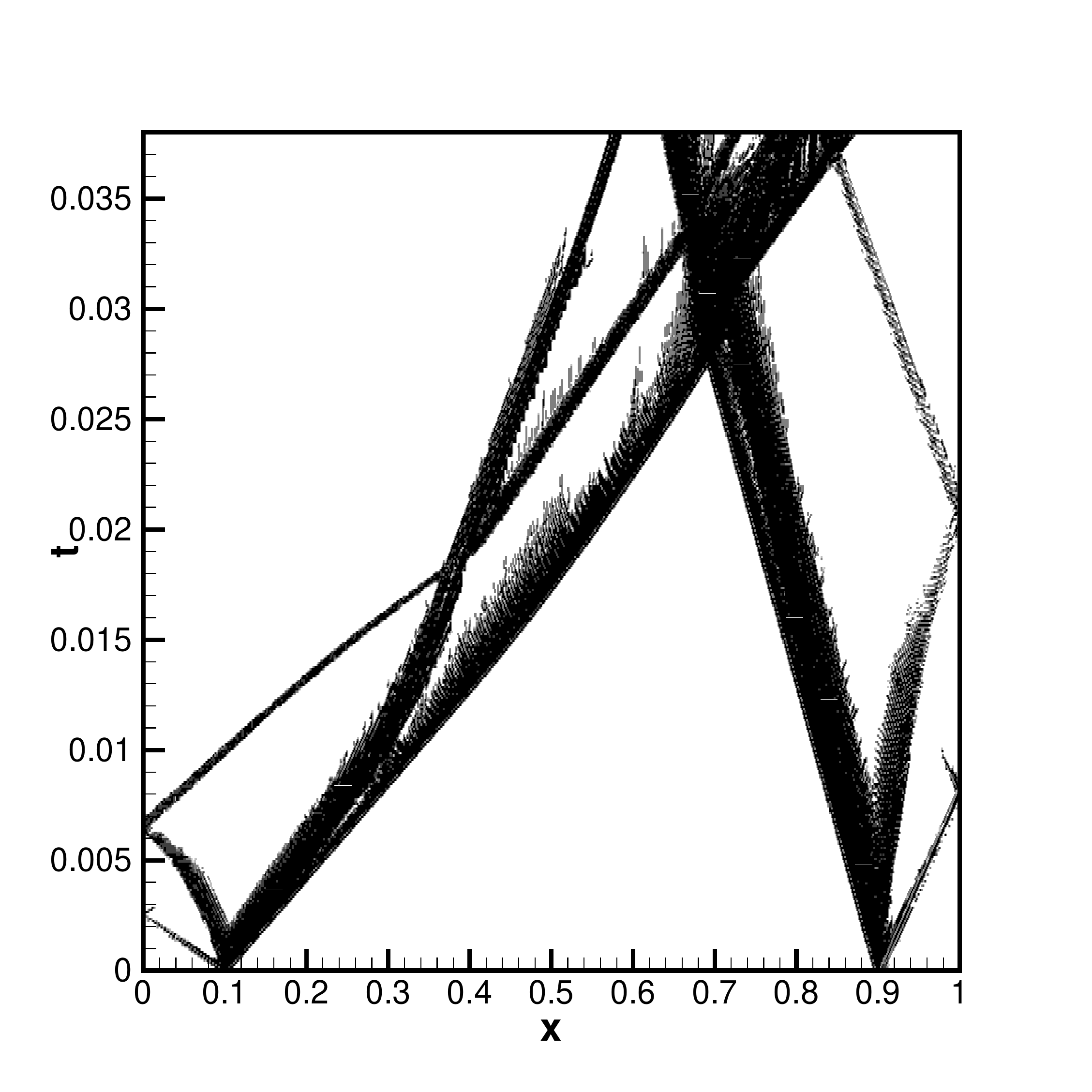}}
  \caption{{
  The blast wave problem.
  $T=0.038$.
  Top:    density, density zoomed in;
  bottom: the locations of the troubled-cells.
  Solid line: the exact solution;
  squares:    the result of the HWENO6      scheme;
  triangles:  the result of the HWENO6-M5-I scheme.
  Number of cells: 800.}
\label{ex7}}
\end{figure}

\noindent {\bf Example 3.7.}
The blast wave problem: one-dimensional Euler equations (\ref{1D-Euler}) with the initial condition:
\begin{equation}
\label{1D-Blast}
(\rho,\mu,p)^T=
\left\{
\begin{array}{c}
 (1,0,10^  3 )^T, \ \   0<x<0.1,\\
 (1,0,10^{-2})^T, \ \ 0.1<x<0.9,\\
 (1,0,10^  2 )^T, \ \ 0.9<x<  1.
\end{array}
\right.
\end{equation}
The computational result of the density $\rho$ and its zoom-in picture obtained by the HWENO6 and HWENO6-M5-I schemes at the final time $T=0.038$ in comparison with the reference ``exact" solution(which is a numerically converged solution computed by the fifth-order finite difference WENO scheme \cite{js1996} with 16000 grid points), as well as the corresponding
locations of the troubled-cells for the HWENO6-M5-I scheme over time are plotted in the Fig \ref{ex7}.
\smallskip

\bigskip

\begin{figure}[htbp]
  \centering{\includegraphics[width=2.5 in,height=2.5 in]{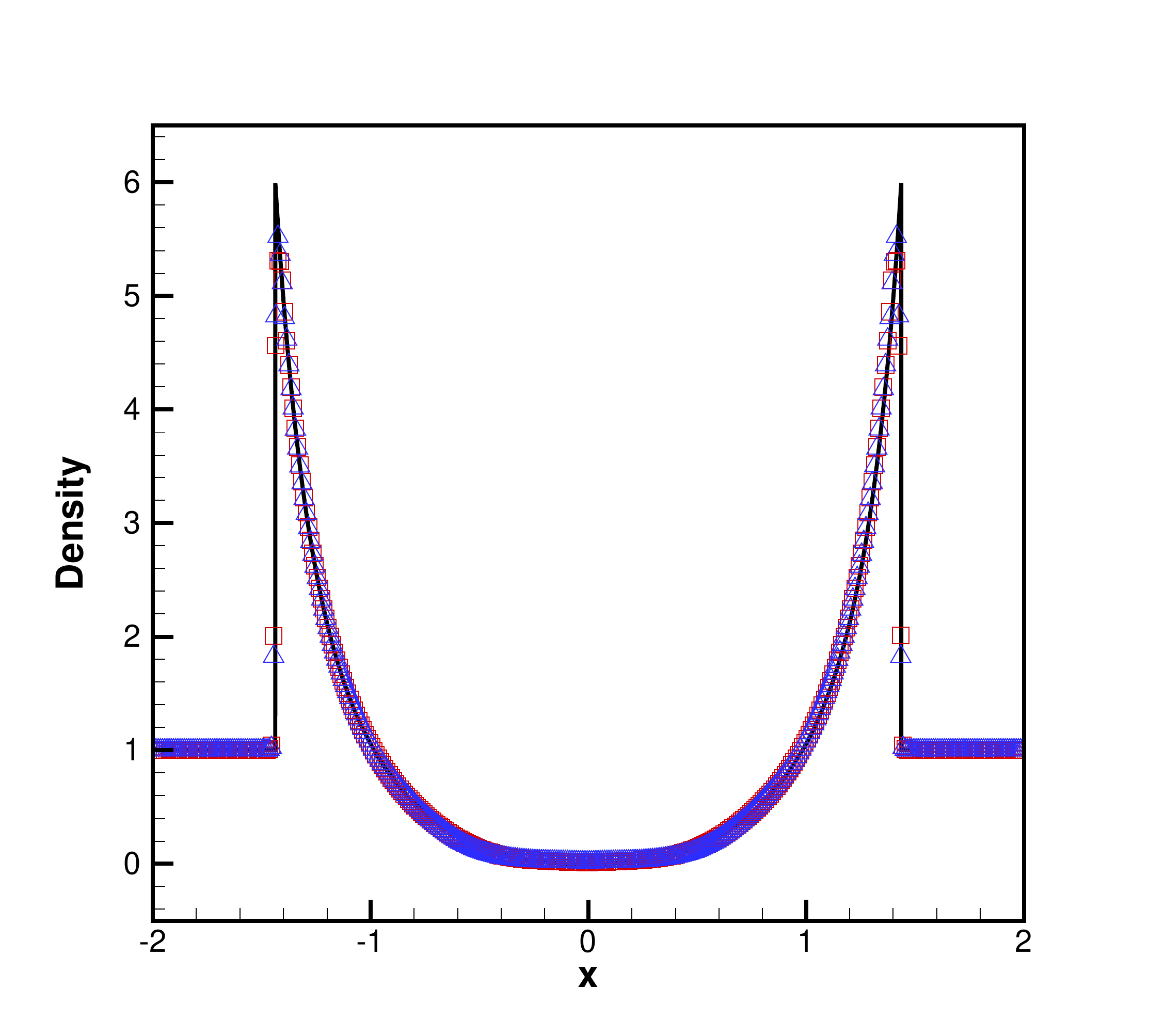}
             \includegraphics[width=2.5 in,height=2.5 in]{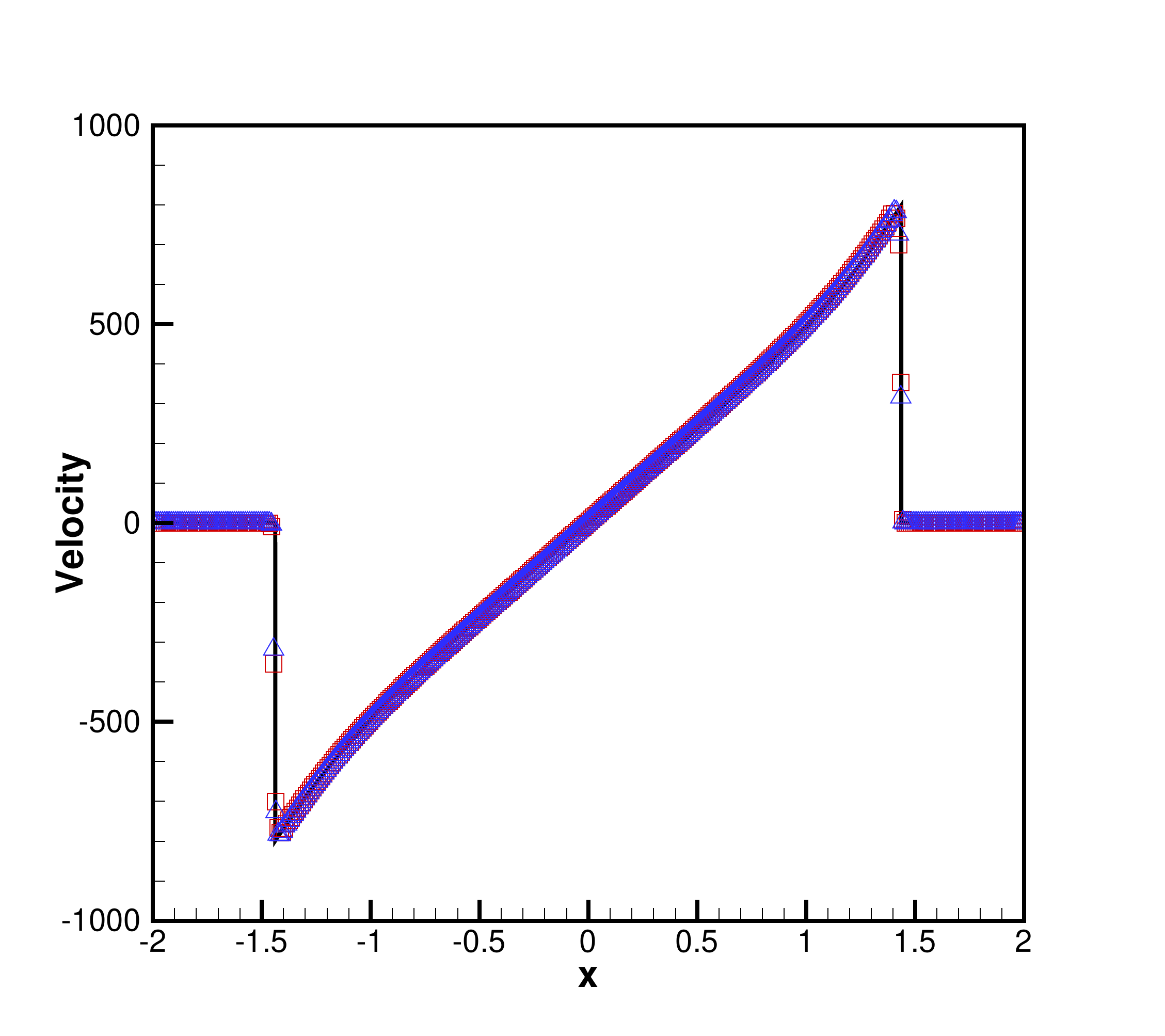}
             \includegraphics[width=2.5 in,height=2.5 in]{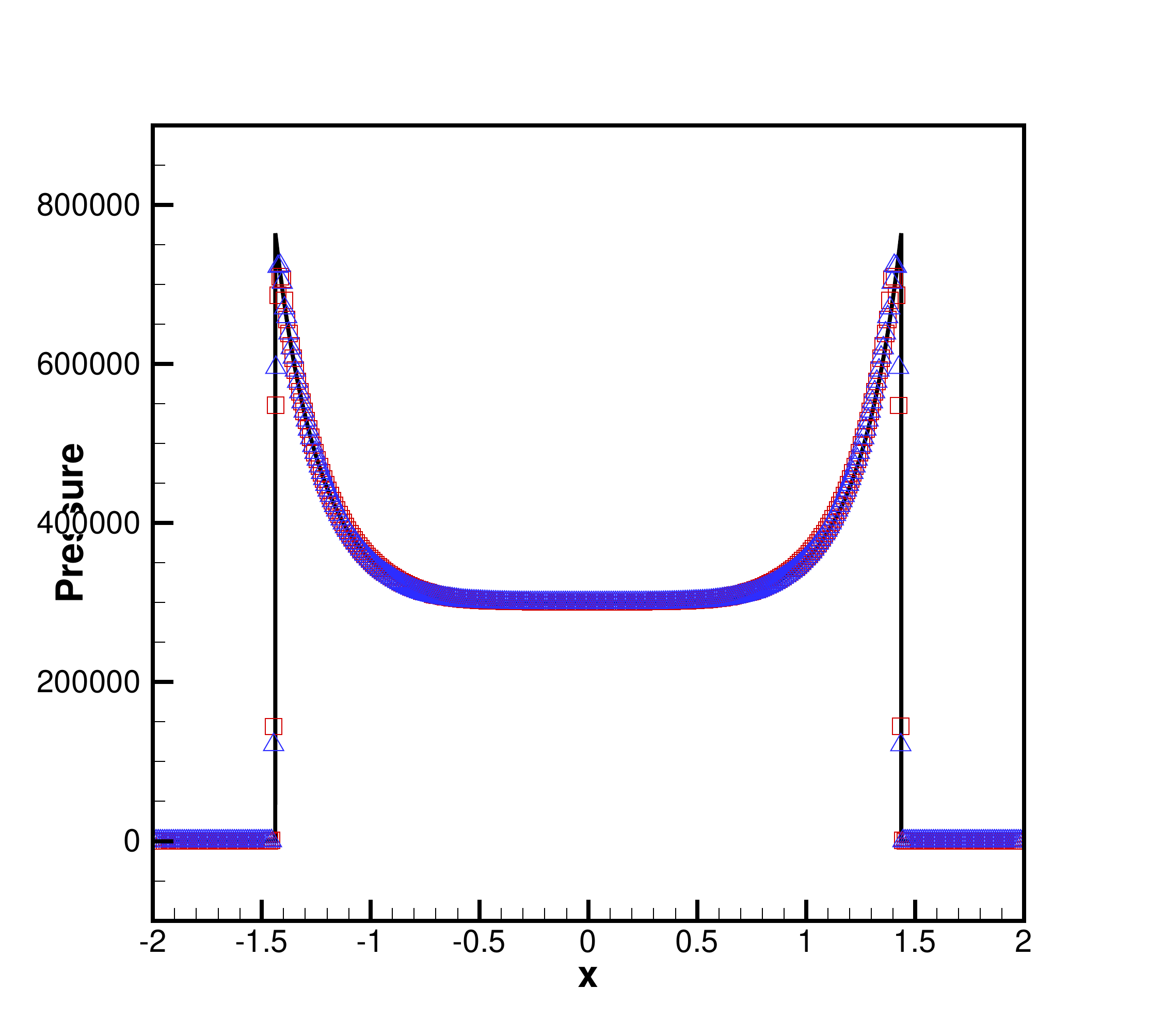}
             \includegraphics[width=2.5 in,height=2.5 in]{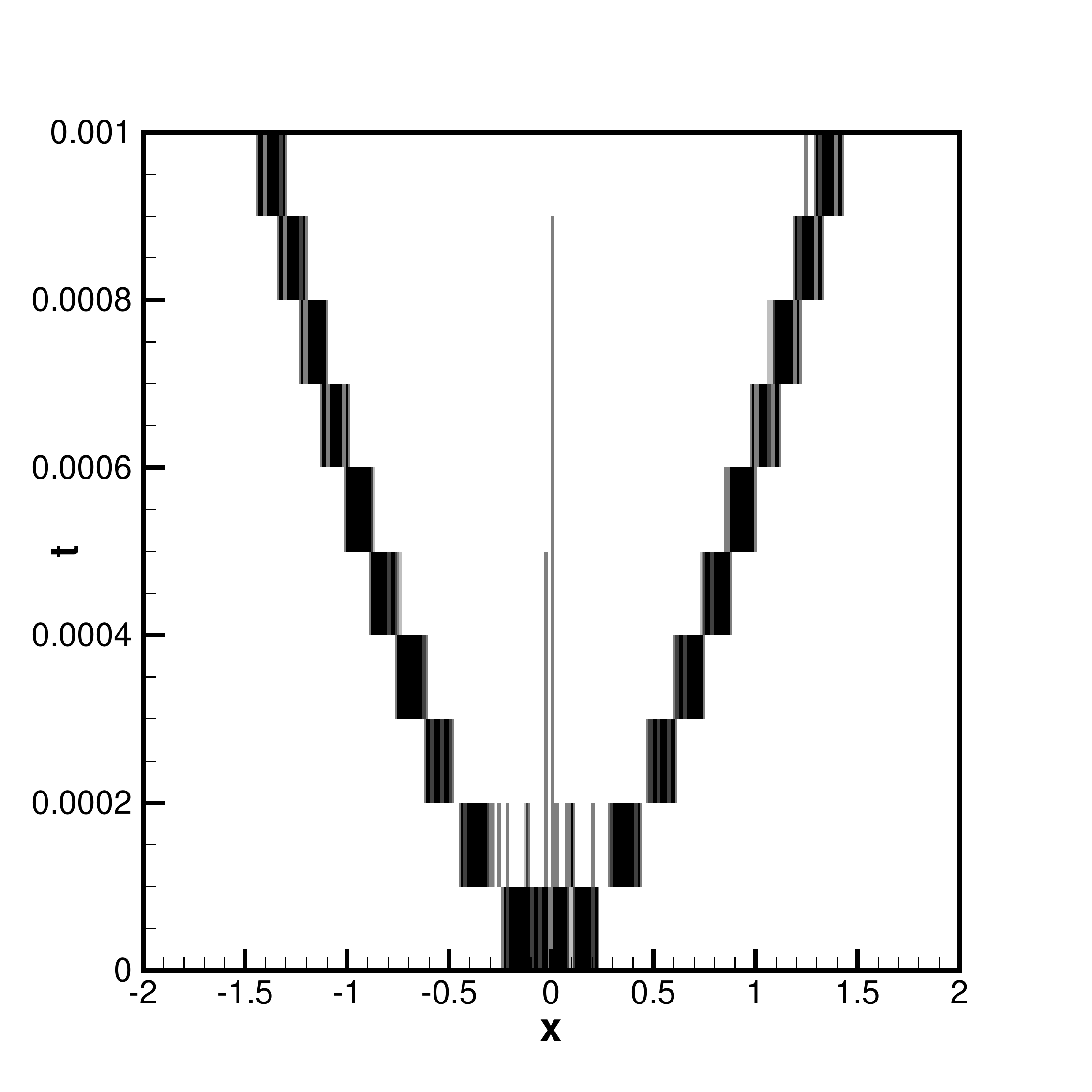}}
  \caption{{
  \footnotesize
  The Sedov blast wave problem.
  $T=0.001$.
  Top:    density, velocity;
  bottom: pressure, the locations of the troubled-cells.
  Solid line: the exact solution;
  squares:    the result of the HWENO6      scheme;
  triangles:  the result of the HWENO6-M5-I scheme.
  Number of cells: 400.}
\label{ex8}}
\end{figure}

\noindent {\bf Example 3.8.}
The Sedov blast wave problem: one-dimensional Euler equations (\ref{1D-Euler}) with the initial condition:
\begin{equation}
\label{1D-Sedov}
(\rho,\mu,E)^T=
\left\{
\begin{array}{c}
 (1,0,10^{-12}                   )^T, \ \
 x \in [-2,2]\setminus the \ \ center \ \ cell,\\
 (1,0,\frac{3200000}{\triangle x})^T, \ \
 x \in the \ \ center \ \ cell.
\end{array}
\right.
\end{equation}
The computed results of the density $\rho$, the velocity $\mu$ and the pressure $p$ obtained by the HWENO6 and HWENO6-M5-I schemes at the final time $T=0.001$ in comparison with the reference exact solution, as well as the corresponding locations of the troubled-cells for the HWENO6-M5-I scheme over time are plotted in the Fig \ref{ex8}.
\smallskip

\bigskip

\begin{figure}[htbp]
  \centering{\includegraphics[width=6.0 in,height=2.0 in]{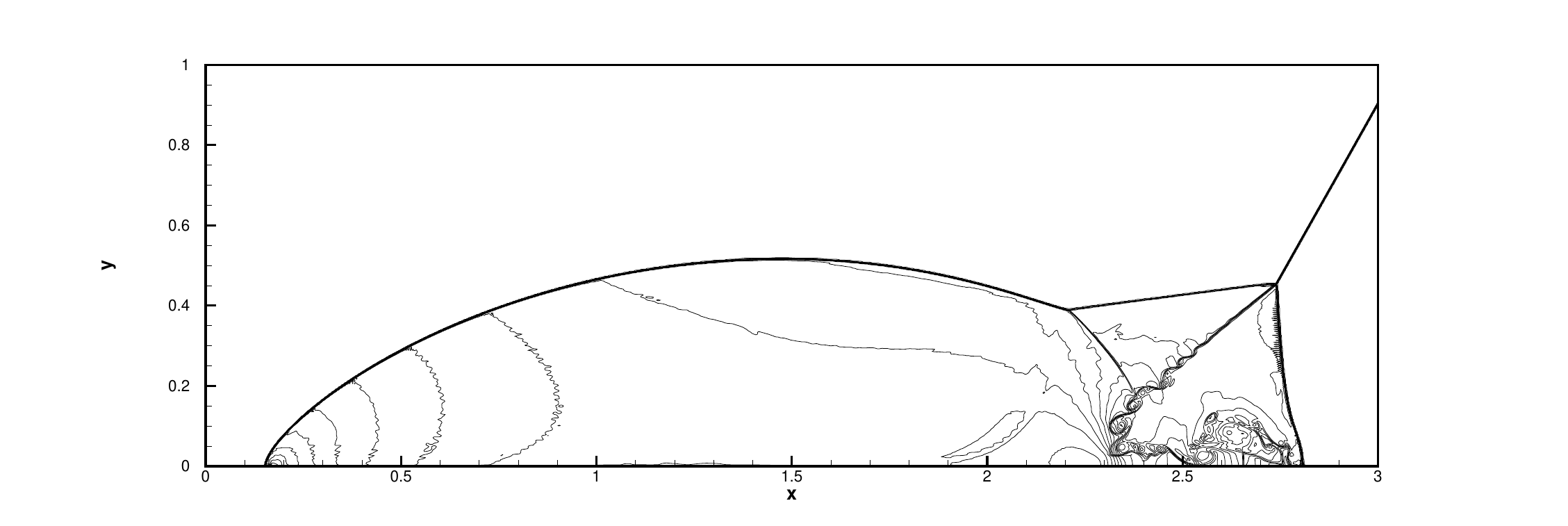}\\
             \includegraphics[width=2.5 in,height=2.5 in]{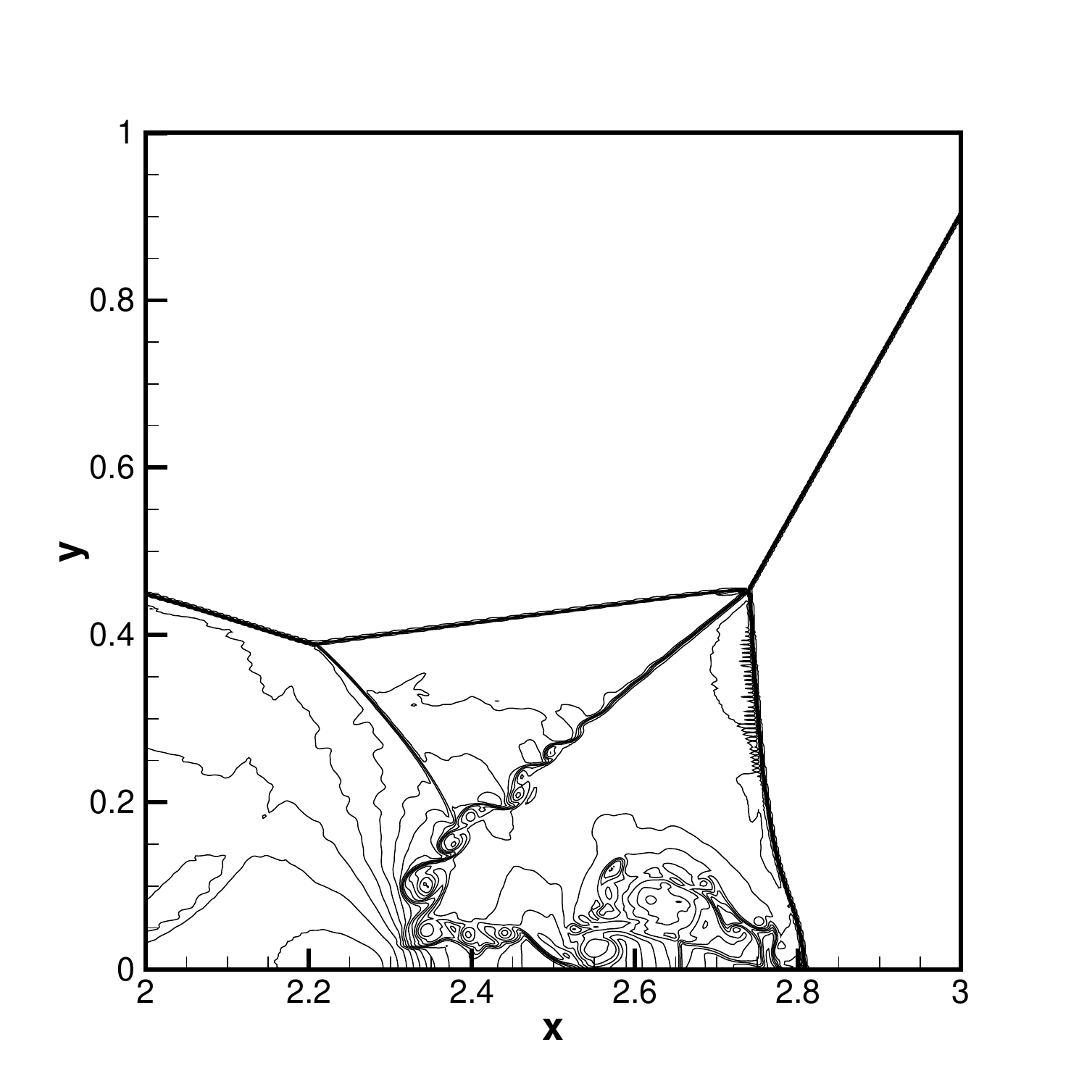}\\
             \includegraphics[width=6.0 in,height=2.0 in]{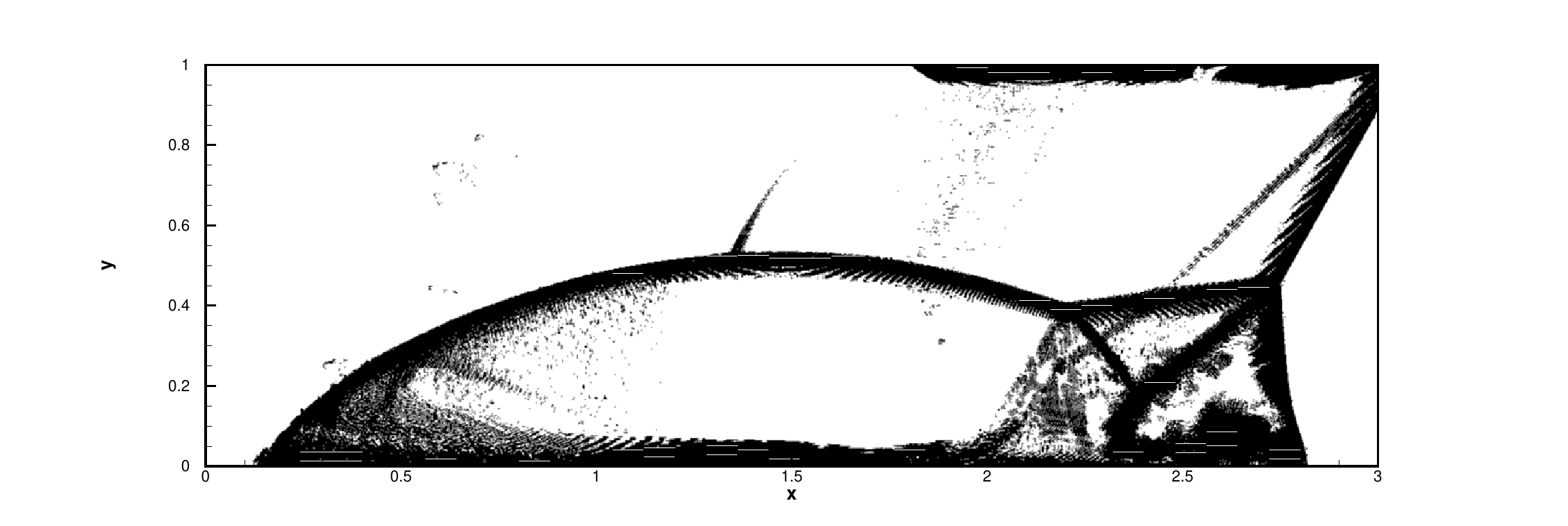}}
  \caption{{
  Double Mach reflection problem.
  $T=0.2$.
  30 equally spaced density contours from 1.5 to 22.7.
  From top to bottom: density contours, zoom-in density contours around the Mach stem, the locations of the troubled-cells at the final time of the HWENO6-M5-I scheme.
  Number of cells: 1600$\times$400 in the region of $[0,4]$$\times$$[0,1]$.}
\label{ex9}}
\end{figure}

\noindent {\bf Example 3.9.}
Double Mach reflection problem: two-dimensional Euler equations (\ref{2D-Euler}) in a computational domain $[0,4]$$\times$$[0,1]$ with a reflection wall lying at the bottom,  starting from the position $(x,y)=(\frac{1}{6},0)$, making a $60^{\circ}$ angle with the $x$-axis.
For the bottom of the domain, at the reflection wall the reflection boundary condition is applied, and at the rest of the bottom the exact post-shock condition is imposed.
For the top    of the domain, the corresponding boundary condition is exactly the  motion of a Mach 10 shock with $\gamma=1.4$.
The contours of the computed density $\rho$ and its blow-up region around the double Mach stem obtained by the HWENO6-M5-I scheme (the reason we do not present the result of HWENO6 scheme is that this scheme for the double Mach reflection problem will 
blow up without modification), as well as the corresponding locations of the troubled-cells for the HWENO6-M5-I scheme at the final time $T=0.2$ are plotted in the Fig \ref{ex9}.
\smallskip

\bigskip

\begin{figure}[htbp]
  \centering{\includegraphics[width=6.0 in,height=2.0 in]{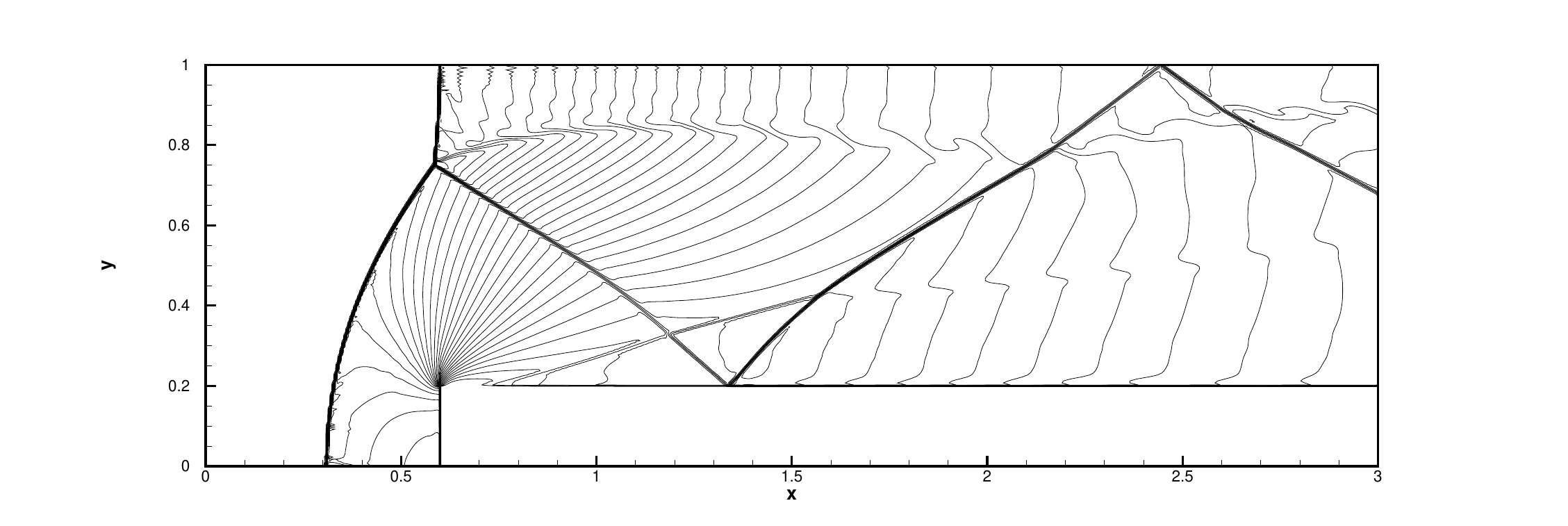}\\
             \includegraphics[width=6.0 in,height=2.0 in]{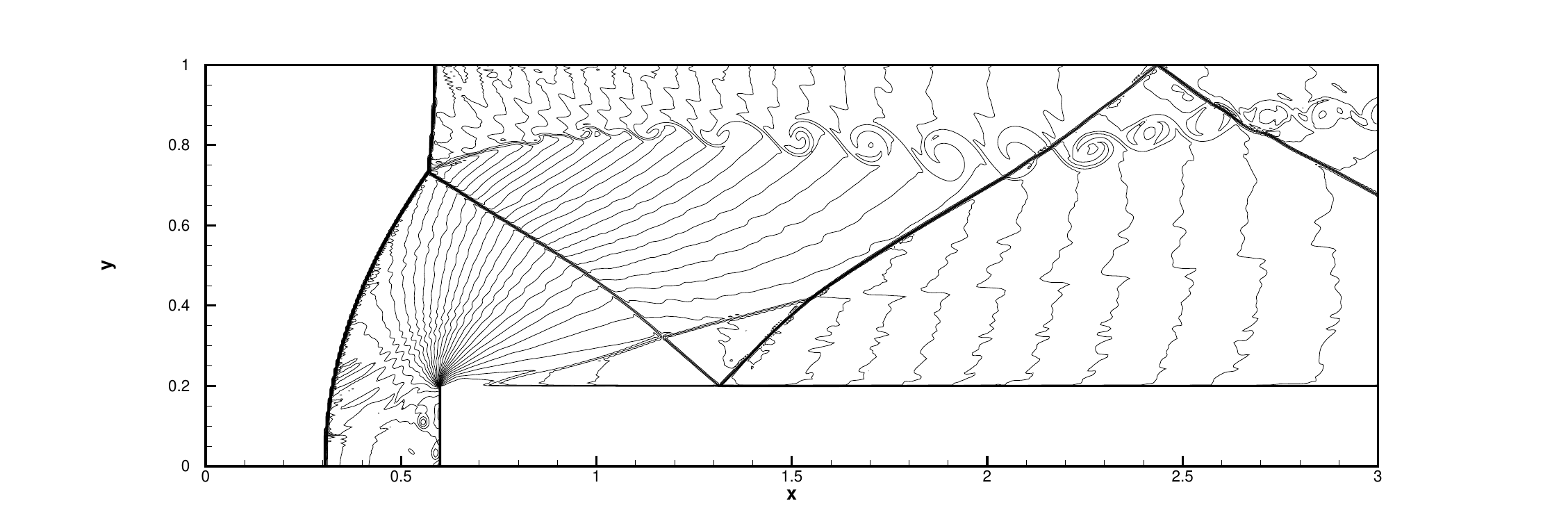}\\
             \includegraphics[width=6.0 in,height=2.0 in]{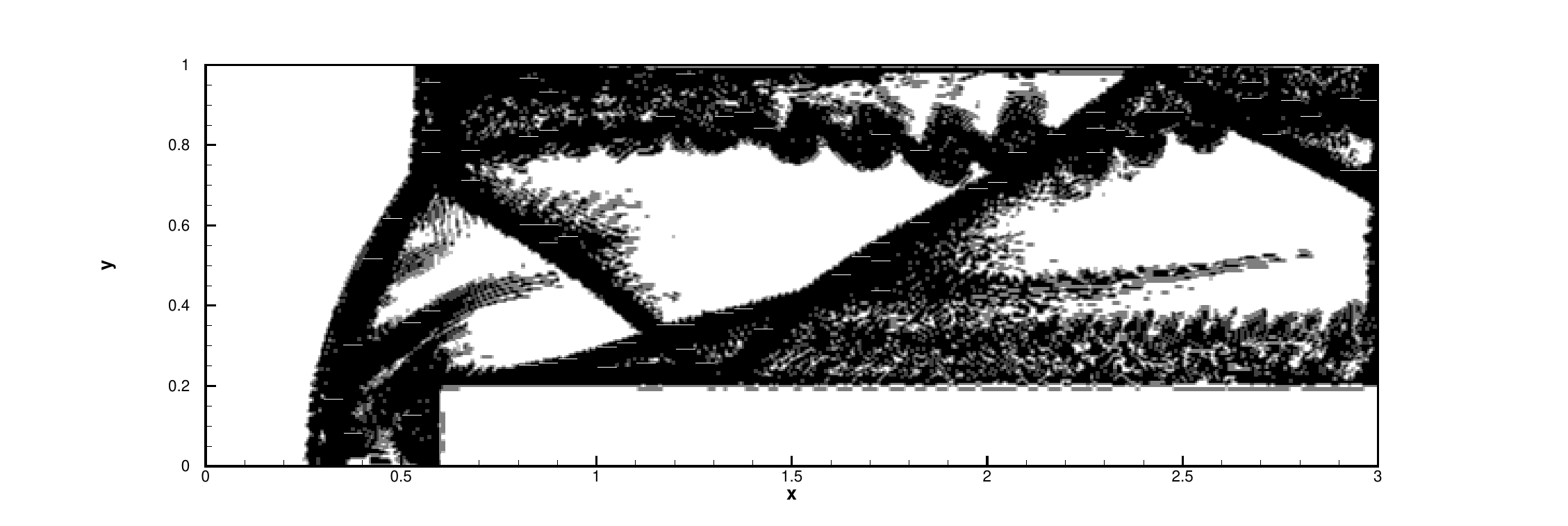}}
  \caption{{
  Forward step problem.
  $T=4$.
  30 equally spaced density contours from 0.32 to 6.15.
  From top to bottom: density contours of the HWENO6 scheme, density contours of the HWENO6-M5-I scheme, the locations of the troubled-cells at the final time of the HWENO6-M5-I scheme.
  Number of cells: 600$\times$200 in the region of $[0,3]$$\times$$[0,1]$.}
\label{ex10}}
\end{figure}

\noindent {\bf Example 3.10.}
Forward step problem: two-dimensional Euler equations (\ref{2D-Euler}) in a one length unit wide and three length units long wind tunnel with a 0.2 length units high step located 0.6 length units from the left side of the tunnel.
At the beginning, we initialize this problem by a right-going Mach 3 flow.
Along the wall of the tunnel the reflection boundary condition is applied, and at the entrance the inflow boundary condition is imposed, while at the exit the outflow boundary condition is imposed.
The contours of the computed density $\rho$ obtained by the HWENO6 and HWENO6-M5-I schemes, as well as the corresponding locations of the troubled-cells for the HWENO6-M5-I scheme at the final time $T=4$ are plotted in the Fig \ref{ex10}.
\smallskip

\bigskip

\noindent {\bf Comment:}
From above six examples, we can see that the results of the HWENO6-M5-I scheme have much better resolutions and sharper shock transitions than those of the HWENO6 scheme without the modification procedure. This might be due to the fact that after modifying the first order moments of the troubled-cells, the proportion of the last two layers is much higher than that of the first two layers in the reconstruction process.
Also, this modification procedure can increase the stability of our scheme according to the double Mach reflection problem.

\section{Concluding remarks}
\label{sec4}
\setcounter{equation}{0}
\setcounter{figure}{0}
\setcounter{table}{0}

In this paper, we have designed a high-order moment-based multi-resolution HWENO scheme for
hyperbolic conservation laws in the one and two dimensional cases on structured meshes.
In comparison with our previous work in \cite{lqs2021}, the new feature of this HWENO scheme is that the zeroth and first order moments rather than the first order derivative are used in the spacial reconstruction algorithm, and only the function values of the Gauss-Lobatto points in one or two dimensional case are needed to be reconstructed.
Also, after the reconstruction algorithm, an extra modification procedure is used to modify those first order moments of the troubled-cells and the corresponding Gauss-Lobatto point values of these troubled-cells need to be updated by repeating the reconstruction algorithm, to enhance both resolution and stability.
At the same time, the linear weights can also be any positive numbers as long as their sum equals one and the CFL number can still be 0.6 for both the one and two dimensional cases.
This HWENO scheme is achieved by reconstructing the Gauss-Lobatto point values, modifying the first order moments of those cells which are identified to be troubled-cells by the KXRCF troubled-cell indicator and repeating the reconstruction algorithm to update the corresponding Gauss-Lobatto point values of these troubled-cells.
In comparison with the multi-resolution WENO scheme, our major advantages are still the
compactness of the stencils and smaller errors under the same meshes when the same order is fixed.
The framework of this moment-based multi-resolution HWENO scheme would be particularly efficient and simple on unstructured meshes, the study of which is
our ongoing work.

{\large \bf Acknowledgments: }The research of C.-W. Shu is partly supported by AFOSR grant FA9550-20-1-0055 and NSF grant DMS-2010107. The research of J. Li and J. Qiu is partly supported by NSFC grant 12071392.  

\end{document}